\magnification=\magstep1
\input amstex
\documentstyle{amsppt}
\catcode`\@=11 \loadmathfont{rsfs}
\def\mycal{\mathfont@\rsfs}
\csname rsfs \endcsname \catcode`\@=\active

\vsize=7.5in

\topmatter 
\title Tight decomposition of  factors \\ and the single generation problem  \\ 
$\text{\it Dedicated to Dan Voiculescu on his 70th birthday}$ \endtitle
\author  Sorin Popa \endauthor

\rightheadtext{Tight decomposition of factors}

\affil     {\it  University of California, Los Angeles} \endaffil

\address Math.Dept., UCLA, Los Angeles, CA 90095-1555\endaddress
\email  popa\@math.ucla.edu\endemail

\thanks Supported in part by NSF Grant DMS-1700344 and the Takesaki Chair in Operator Algebras  \endthanks

\abstract  A  II$_1$ factor $M$ has the {\it stable single generation} ({\it SSG}) property if any amplification  $M^t$, $t>0$, can be generated 
as a von Neumann algebra by a single element. We discuss a conjecture stating that 
if $M$ is SSG, then $M$ has a {\it tight} decomposition, i.e., there exists a pair of hyperfinite II$_1$ subfactors $R_0, R_1  \subset M$ 
such that $R_0 \vee R_1^{op}=\Cal B(L^2M)$.   We provide supporting evidence, explain why the conjecture is 
interesting, and discuss possible approaches to  settle it.   We also prove some related results.

\endabstract 

\endtopmatter

\document

\heading 1. Introduction \endheading

It has been recently conjectured in (5.1(b) in [P18]; see also Section 7 in [P19])
that if a II$_1$ factor $M$ is {\it stably single generated} (abbreviated hereafter as {\it SSG}), i.e., if $M^t$ is single 
generated as a von Neumann algebra for any $t>0$, then $M$ has an $R$-{\it tight decomposition}, meaning that  
it contains  hyperfinite subfactors $R_0, R_1 \subset M$ such that $R_0 \vee R_1^{op}=\Cal B(L^2M)$. A weaker conjecture in (5.1(a) of [P18]) states that if $M$ is 
SSG then it admits a {\it properly infinite} $R$-{\it pair},  i.e., hyperfinite  subfactors $R_0, R_1\subset M$ so that $R_0\vee R_1^{op}$ is a properly  infinite von Neumann algebra in $\Cal B(L^2M)$. 

One can easily see that if $M$ has non-trivial fundamental group, $\Cal F(M)\neq 1$, then $M$ is SSG if and only if it is finitely generated. Thus, 
since the II$_1$ factor $L(\Bbb F_\infty)$ of the free group with infinitely many generators $\Bbb F_\infty$ has non-trivial fundamental group by [V88] 
(in fact one even has $\Cal F(L(\Bbb F_\infty))=\Bbb R_+$ by [R91]), if $L(\Bbb F_\infty)$ would be finitely generated 
and the implication ``SSG $\Rightarrow$ $R$-tight'' holds true, then $L(\Bbb F_\infty)$ would follow $R$-tight, 
contradicting a result in ([GP98]). More precisely, the result in [GP98], whose proof is  based on 
Voiculescu's major breakthrough techniques and free entropy theory in ([V96]), shows that if $M=L(\Bbb F_n)$, $2\leq n \leq \infty$, 
then $M$ is not {\it weakly thin}, meaning that one cannot find two AFD subalgebras $B_0, B_1\subset M$ and $X\subset L^2M$ finite such that  
sp$B_0XB_1$ is dense in $L^2M$. In particular, there exists no pair of hyperfinite factors $R_0, R_1\subset M$ so that $R_0\vee R_1^{op}\subset \Cal B(L^2M)$ admits a 
finite cyclic set. But if $R_0\vee R_1^{op}=\Cal B(L^2M)$ then any non-zero vector in $L^2M$ is cyclic. Even if  
$R_0 \vee R_1^{op}$ is merely properly infinite, it always has a cyclic vector (see e.g., [D57], [S71]). 

So if true, these conjectures would  imply that $L(\Bbb F_\infty)$ cannot be generated by finitely many elements (which is what one calls being {\it infinitely generated}). By (Corollary 4.7 in [R92]), this in turn   
would imply that the free group factors $L(\Bbb F_n), 2 \leq n \leq \infty$, are all non-isomorphic. 

The above conjectures have been triggered by a result in [P18], showing that any separable II$_1$ factor $M$ has a {\it coarse decomposition}, 
in the sense that there exist embeddings of the hyperfinite II$_1$ factor, $R_0, R_1\hookrightarrow  M$, such that $R_0\vee R_1^{op}$ is 
finite (thus isomorphic to $R_0\overline{\otimes} R_1^{op}$). More precisely, they were motivated 
by the method we used to prove this result: the  coarse pair of hyperfinite II$_1$ subfactors $R_0, R_1$ in $M$ is constructed 
recursively, as inductive limits of dyadic finite dimensional factors $R_{0,n}\nearrow R_0$, $R_{1,n}\nearrow R_1$, so that at each step $n$ more and more of the vectors in a 
countable dense subset $\Cal L \subset L^2M$ implement asymptotically  a specific type of state on $R_{0,n}\vee R_{1,n}^{op} \simeq R_{0,n}\otimes R_{1,n}^{op}$, namely the trace $\tau\otimes \tau$. 

It is reasonable to believe that one can make this ``iterative construction with constraints''  so that all the vectors in $\Cal L \subset L^2M$ 
implement  asymptotically  states that   ``stay away'' 
from $\tau\otimes \tau$.  However, for certain factors this is not possible: one can easily deduce from (Theorem 4.2 in [GP98]) that if $M$ is a free group factor 
then any choice of an increasing sequence of  dyadic factors $R_{0,n}, R_{1,n}$ ends up producing 
a pair of hyperfinite factors $R_0, R_1\subset M$ with $R_0\vee R_1^{op}$  having a coarse part  (see Theorem 2.9 below). 
In other words, no matter what one does, some of the vectors in $L^2M$ will necessarily  implement $\tau\otimes \tau$ on $R_0\vee R_1^{op}$,  
making it impossible for  $R_0\vee R_1^{op}$ to be properly infinite. 

To escape this ``coarseness trap'', the iterative construction should thus take into account  ``special properties'' that $M$ may have. 
The above conjectures come from our belief that the SSG property of $M$ is enough to insure that one can build iteratively a pair of hyperfinite factors $R_0, R_1\subset M$ 
so that $R_0\vee R_1^{op}$ is properly infinite, in fact even tight. 

Our purpose in this paper is to discuss in more details these conjectures and the properties involved, especially tightness and SSG, and prove some related results. 

We begin by introducing in Section 2 
some terminology about bimodule decompositions of II$_1$ factors $M$ over pairs of hyperfinite subfactors $R_0, R_1\subset M$ (respectively over a hyperfinite 
subfactor $R\subset M$), 
according to the type and other properties  of the von Neumann algebra $R_0 \vee R_1^{op}\subset \Cal B(L^2M)$ 
(respectively the type of $R\vee R^{op}$ on $L^2(M\ominus R)$): finite/coarse, with its  corresponding multiplicity (coupling constant); factorial but  
properly infinite, with its subclasses by type (I, II$_\infty$, III$_\lambda$, $0\leq \lambda \leq 1$).  We provide examples and 
prove some basic criteria. We formulate problems, underlying the idea that the study of $R$-bimodule decomposition properties 
of a II$_1$ factor is both interesting and non-trivial,  that's worth investigating in its own right. 
Then in Section 3 we concentrate on tight decompositions of a II$_1$ factor $M$, proving some criteria and giving examples. 

In Sections 4 and 5 we discuss the SSG property for 
II$_1$ factors. We notice that SSG has good permanence properties, being preserved under crossed products, quasi-regular inclusions, finite index restrictions/extensions, 
inductive limits, etc. These properties are immediate consequences of results in ([Sh05], [DSSW07]), but we give a different, self-contained treatment, for the reader's convenience. 
In particular, we provide a short argument to a result in ([DSSW07]), showing that if one denotes 
by $\text{\rm ng}(M)\in \{2,3,...\}\cup \{\infty\}$ the minimal number of self-adjoint elements that can generate the II$_1$ factor $M$, 
then either $\text{\rm ng}(M^t)=\infty$ for all $t>0$, 
or $\text{\rm ng}(M^t)= O(t^{-2}+1)$ (note that this shows  that if $\underset{t\rightarrow 0}\to{\liminf} \ t^2\text{\rm ng}(M^t)=0$ then $M$ is SSG). 
More precisely, we use Voiculescu's ``efficient counting'' of generators for $M^{1/k}$ from the ones of $M$ in ([V88]), 
to show that ng$(M)\leq n$  if and only if ng$(M^{1/k})\leq  (n-1)k^2 +1$.  

In the last Section 6 we discuss various possible strategies for constructing  $R$-tight decompositions, or more generally 
properly infinite $R$-decompositions, from the SSG property. We also formulate a number of problems.   

\vskip.05in 

For general notations and basic facts about  II$_1$ factors that we use in this paper, we refer the reader to ([AP17], [P18]). 
As usual, we use the notation $R$ to designate the hyperfinite II$_1$ factor $R=(\Bbb M_2(\Bbb C), tr)^{\overline{\otimes}\infty}$, which 
is the unique (separable) AFD II$_1$ factor ([MvN43]), in fact even the unique amenable II$_1$ factor ([C76]).  But we will mostly give this notation a generic meaning, as the      
adjective ``hyperfinite'' (for instance an $R$-pair, means a pair of hypefinite II$_1$ factors).

\heading 2.  Bimodule decomposition of II$_1$ factors 
\endheading

We fix here some terminology about the decomposition of a given II$_1$ factor $M$ as a Hilbert bimodule over its subfactors. 
We are interested in two such cases: viewing $M$ as a bimodule $_QL^2M_Q$ over a single subfactor $Q\subset M$,  and as a bimodule $_QL^2M_P$ 
over a pair of subfactors $Q, P\subset M$. 
We will distinguish these bimodules by properties (such as type) of the von Neumann algebra $Q\vee Q^{op}$, respectively 
$Q\vee P^{op}$, generated in $\Cal B(L^2M)$ by the operators of left and right multiplication by $Q$, respectively left multiplication by $Q$ and right multiplication 
by $P$. Since we always have $_QL^2M_Q=L^2Q \oplus L^2(M\ominus Q)$, with the projection $e_Q$ lying in the center of $Q\vee Q^{op}$ and $Q\vee Q^{op}e_Q=\Cal B(L^2Q)$, 
these properties  (such as coarseness) will actually refer to the ``interesting part'' $_QL^2(M\ominus Q)_Q$. 

We will denote by $(Q\vee Q^{op})_{fin}, (Q\vee Q^{op})_\infty$ (resp. $(Q\vee P^{op})_{fin}, (Q\vee P^{op})_\infty$) 
the finite and properly infinite direct summands of $Q\vee Q^{op}$ (resp. $Q\vee P^{op}$). 

We are particularly interested in the case when $Q\vee Q^{op}$ and $Q\vee P^{op}$ are {\it homogeneous} von Neumann algebras (the former on $L^2(M\ominus Q)$), 
of either {\it finite} of {\it properly infinite} type. One should note that since $Q, P$  are II$_1$ factors, the finite part $(Q\vee Q^{op})_{fin}$, 
resp. $(Q\vee P^{op})_{fin }$, is in fact a II$_1$ factor of the form $Q\overline{\otimes} Q^{op}$, resp. $Q\overline{\otimes} P^{op}$. So finite homogeneous means 
$Q\vee Q^{op}(1-e_Q) \simeq Q\overline{\otimes} Q^{op}$, resp. $Q\vee P^{op}\simeq Q\overline{\otimes} P^{op}$,  which in the terminology of ([P18]) amounts 
to $Q$ being {\it coarse} in $M$, resp. $Q, P$ being a {\it coarse pair} in $M$. If this is the case,  
then we call the coupling constant $\text{\rm dim}_{Q\overline{\otimes} Q^{op}}L^2(M\ominus Q)$ 
of $Q\vee Q^{op} \subset \Cal B(L^2(M\ominus Q))$, resp. $\text{\rm dim}_{Q\overline{\otimes} P^{op}}L^2M$ of  $Q \vee P^{op}$  in $\Cal B(L^2M)$,  the 
{\it multiplicity} of the coarse embedding $Q\subset M$, resp. of the coarse pair $Q, P\subset M$, and denote it $c(M;Q, Q)$, resp $c(M;Q, P)$. So apriori, this is a number in $[0,\infty]$.  

If $Q\subset M$ (resp. $Q, P \subset M$) is coarse with $Q$ (resp. $Q, P$) hyperfinite,  then we call it an $R$-{\it coarse subfactor} 
(resp. an $R$-{\it coarse pair}), or that it gives an $R$-{\it coarse decomposition of $M$} (resp. {\it $R$-coarse pair decomposition} of $M$). 

The $Q$-bimodule decomposition  (resp. $Q-P$ bimodule decomposition) 
of $M$ is {\it properly infinite} if $Q\vee Q^{op}\subset \Cal B(L^2(M\ominus Q))$, resp $Q \vee P^{op}\subset \Cal B(L^2M)$ is a 
{\it properly infinite} von Neumann algebra. Note that this is equivalent to these bimodules having 
no coarse part. Such a properly infinite bimodule decomposition is of {\it type} I, II$_\infty$, III, if $Q\vee Q^{op}$, resp. $Q\vee P^{op}$, is of type I, II$_\infty$,  III, 
adding {\it factorial} when this algebra is a factor. In case $Q\vee Q^{op}$ (resp. $Q\vee P^{op}$) is a factor of type III$_\lambda$, $0\leq \lambda \leq 1$, we say that $Q\subset M$ (resp. $Q, P \subset M$) 
gives a III$_\lambda$-{\it factorial decomposition} of $M$. 

We also single out the case when $Q\vee Q^{op}\subset \Cal B(L^2(M\ominus Q))$, resp. 
$Q \vee P^{op}\subset \Cal B(L^2M)$,  is {\it cyclic}, meaning that it has a cyclic vector. Note that if $Q, P \simeq R$ and the cyclic vector is $\hat{1}\in L^2M$, 
then the terminology used in ([P94], [P97], [GP98]) is that $M$ is a {\it thin} factor. 

Note that a pair $Q, P \subset M$ is cyclic if and only if the finite part of $Q\vee P^{op}\subset \Cal B(L^2M)$ is cyclic, 
a condition that's equivalent to having multiplicity (or coupling constant)  $\leq 1$. 
In particular, if $Q, P$ give a properly infinite decomposition of $M$, then it is automatically cyclic (see e.g., [S71]). However, as we will see in  Example 2.2.4$^\circ$ and Corollary 2.3,   
a II$_1$ factor $M$ may have a cyclic $R$-coarse pair decomposition, but nevertheless be non-amenable, even non-Gamma.  

A pair of hyperfinite subfactors $R_0, R_1\subset M$ gives a {\it weakly $R$-thin decomposition} of $M$ if there exists a finite set $X\subset L^2M$ 
such that $[R_0XR_1]=L^2M$. Note that this condition is equivalent to the fact that the coupling constant of $(R_0 \vee R_1^{op})_{fin}$ is finite 
(namely majorized by $|X|$). 

The next result shows that if $R_0, R_1\subset M$ is an $R$-bimodule decomposition of $M$, then all the characteristics of this decomposition (like type, 
multiplicity) do not depend on the individual unitary conjugacy class of $R_0$ and $R_1$, nor on taking amplifications $R_0^t, R_1^t \subset M^t$ 
by $t>0$.  The meaning here of the amplification by a non integer $t>0$ of an inclusion of II$_1$ factors is, as usual, viewed modulo unitary conjugacy 
of the subfactor.

\proclaim{2.1. Lemma} Let $M$ be a $\text{\rm II}_1$ factor and $R_0, R_1\subset M$ a pair of hyperfinite subfactors. 
Let $R_0\vee R^{op}_1=\Cal R_{fin}\oplus \Cal R_\infty \subset \Cal B(L^2M)$, with $\Cal R_{fin}=(R_0\vee R_1^{op})_{fin}  \simeq R_0\overline{\otimes} R_1^{op}$ and $\Cal R_\infty=(R_0\vee R_1^{op})_\infty$. 

\vskip.05in
$1^\circ$ If $u, v \in \Cal U(M)$ and one denotes $P_0=uR_0u^*, P_1=vR_1v^*$, then the decomposition $P_0\vee P_1^{op}=\Cal P_{fin}\oplus \Cal P_\infty$, 
where $\Cal P_{fin}=(P_0\vee P_1^{op})_{fin}$, $\Cal P_\infty=(P_0\vee P_1^{op})_\infty$, 
is spatially conjugate to $\Cal R_{fin}\oplus \Cal R_\infty$ via $\text{\rm Ad}(uv^{op})$ inside $\Cal B(L^2M)$. 

\vskip.03in

$2^\circ$ If $p\in \Cal P(R_0)\cap \Cal P(R_1)$ and one denotes $M^p=pMp$, $R_0^p=pR_0p$, $R_1^p=pR_1p$, then 
 with the identifications  $\Cal B(L^2(M^p))=\Cal B(L^2(pMp))=pp^{op}\Cal B(L^2M)pp^{op}$, the decomposition into finite and properly 
infinite parts $R_0^p\vee {R_1^p}^{op} =\Cal R^p_{fin}\oplus \Cal R^p_\infty \subset \Cal B(L^2(M^p))$ satisfies $\Cal R^p_{fin}=pp^{op}\Cal R_{fin}pp^{op}$, 
$\Cal R_\infty^p=pp^{op}\Cal R_\infty pp^{op}$. Moreover, the coupling constant of $\Cal R_{fin}$ is equal to the coupling constant of $\Cal R_{fin}^p$. 

\endproclaim
\noindent
{\it Proof}. Part $1^\circ$ is trivial. To prove part $2^\circ$, let $\xi_1, ..., \xi_n\in L^2M$ be so that $[R_0\xi_iR_1]$ are mutually 
orthogonal, with $\xi_1, \xi_2, ..., \xi_{n-1}$ implementing $\tilde{\tau}=\tau\otimes \tau$ on $R_0\vee R_1^{op}$ and $\xi_n$ implementing 
$\tilde{\tau}( \cdot \ q)$ for some projection $0\neq q\in R_0\subset R_0\overline{\otimes} R^{op}_1=(R_0\vee R_1^{op})q_{fin}$, 
where $q_{fin}=\sum_i [R_0\xi_i R_1]$ gives the support of the finite part of $R_0\vee R_1^{op}$. Note that  with these notations, the coupling constant 
of $(R_0\vee R_1^{op})_{fin}$ is equal to $n-1+\tau(q)$. Note also that due to the factoriality of the finite part $\Cal R_{fin} \simeq R_0\overline{\otimes} R_1^{op}$ 
(which we already used to take $q\in R_0$), one can assume $q\in R_0$ commutes with $p$ and satisfies $\tau(pq)=\tau(p)\tau(q)$. 

Note that the  support $q^p_{fin}$ of the finite part of $R^p_0\vee {R_1^p}^{op}$ is given by $q^p_{fin}=pp^{op}q_{fin}pp^{op}$. Moreover, if we denote 
$\eta_i=p\xi_ip \in L^2(M^p)$, $1\leq i \leq n$,  then $q^p_{fin}=\sum_i [R_0^p\eta_i R_1^p]$.  
Also, for each $1\leq i \leq n-1$, $\eta_i$ implements the trace $\tilde{\tau}_p=\tau_p \otimes \tau_p$ on $\Cal R^p_{fin}=R_0^p\overline{\otimes} {R_1^p}^{op}$  
while $\eta_n$ implements $\tilde{\tau}_p(\cdot \ q)$. This shows that the coupling constant of $\Cal R_{fin}^p$ is equal to $n-1 + \tau_p(q)=n-1 + \tau(q)$. 
\hfill 
$\square$ 

\proclaim{2.2. Corollary} Let $M$ be a $\text{\rm II}_1$ factor and $R_0, R_1 \subset M$ a pair of hyperfinite $\text{\rm II}_1$ subfactors. 
Then the type of the corresponding $R_0-R_1$ bimodule decomposition of $M$ and the coupling constant of 
$(R_0 \vee R_1^{op})_{fin}$ do not depend on the unitary conjugacy classes of $R_0, R_1$, nor on taking amplifications $(R_0, R_1 \subset M)^t$, $t>0$. 

In particular, if $R_0, R_1$ give a coarse, respectively properly infinite, decomposition of $M$, then $R_0^t, R_1^t$ give a 
coarse, respectively properly infinite, decomposition of $M^t$, $\forall t>0$, which in the coarse case has the same multiplicity, $\forall t>0$. 
\endproclaim 
\noindent
{\it Proof}. If $0< t \leq 1$, then all statements are direct consequences of  Lemma 2.1. Moreover, by applying Lemma 2.1 to $\Bbb M_k(R_0), \Bbb M_k(R_1)\subset \Bbb M_k(M)$ 
and $t=1/k$, one obtains that the statement holds true for $t=k$ an integer as well. Combining with the case $t\leq 1$,  it follows that it actually holds true $\forall t>0$. 

\hfill 
$\square$

\vskip.1in

\vskip.05in 
\noindent
{\bf 2.3. Examples}  $1^\circ$ if $Q$ is a II$_1$ factor, $\Gamma \curvearrowright Q$ is a free action of a countable group and $M=Q \rtimes \Gamma$, 
then $Q\vee Q^{op}\subset \Cal B(L^2(M\ominus Q))$ of type I$_\infty$ with atomic center, and cyclic.  

\vskip.05in 

$2^\circ$ If $M=Q \overline{\otimes} P$, then $Q=Q\otimes 1, P=1\otimes P$ is a coarse, cyclic pair in $M$. If $Q_0 \subset Q, P_0\subset P$ are II$_1$ subfactors, 
then $Q_0, P_0$ is also a coarse pair in $M$, and its multiplicity  is equal to $[Q:Q_0] [P:P_0]$. So by [J82], the multiplicity of a coarse pair can take any value in the set $\{4\cos^2 \pi/n \mid n \geq 3\}\cup [4, \infty]$. 
But we have no example where the multiplicity is less than $1$, or in the remaining range above $1$.

\vskip.05in 

$3^\circ$ If $M$ is a non-prime II$_1$ factor of the form $Q_1 \overline{\otimes} P_1$ and $P_0\subset P_1$ is an irreducible subfactor, 
then $Q=Q_1 \overline{\otimes} P_0$, $P=1\otimes P_1$ is a II$_\infty$ factorial pair 
(thus also cyclic). 

\vskip.05in 

$4^\circ$ Let $\Gamma$ be a group with subgroups $H, \Gamma_0\subset \Gamma$ satisfying 
$H\Gamma_0=\Gamma$. Let $(B, \tau)$ be a tracial von Neumann algebra and $\Gamma \curvearrowright^\sigma B$ a trace preserving action. 
Let also $B_0\subset B$ be $\Gamma_0$-invariant and $B_1\subset B$ be $H$-invariant subalgebras 
such that $B=\overline{\text{\rm sp}B_1B_0}$. Assume that $M= B \rtimes \Gamma$, $B_1 \rtimes H$, $B_0 \rtimes  \Gamma_0$ are all factors (which happens 
if either the group involved is ICC and the action is ergodic on the center of the algebra it acts on, 
or if the action is free ergodic). Then $M=B \rtimes \Gamma$ is a $\text{\rm II}_1$ factor and $Q=B_1\rtimes H, P=B_0\rtimes \Gamma_0$ is a cyclic pair of subfactors. 

If in addition $H \cap \Gamma_0=\{e\}$ and $B_0 \perp B_1$,  then the pair $(Q, P)$ is coarse. 

Note that by (Remark 1.8 in [GP98]), one can take 
$\Gamma$ non-amenable with $H, \Gamma_0$ amenable and $\Gamma_0 \cap H=\{e\}$. More precisely, if one takes $k$ be the field of all algebraic real numbers 
and $\Gamma=PSL(2, k)$, then $\Gamma$ is ICC non-inner amenable which contains the subgroups $H=SO(2, k)=\{u\in \Gamma \mid u^tu=1\}$, 
$\Gamma_0$ the upper triangular matrices in $\Gamma$, and we have indeed $\Gamma = H \Gamma_0$, $\Gamma_0 \cap H=\{1\}$, with $H$ abelian, $\Gamma_0$ 
amenable ICC (see [B05] for a detailed discussion and more examples of groups having such a ``tight decomposition'', 
known as {\it Zappa-Szep product} decomposition\footnote{I am grateful to the referee  
for pointing this out to me}). Thus, if one takes 
$\Gamma \curvearrowright B=\overline{\text{\rm sp} B_0B_1}$ to be any free action on an AFD algebra which has the above properties,  with $B_0 \perp B_1$, 
then $Q=B_1\rtimes H$, $P=B_0 \rtimes \Gamma_0$ gives a hyperfinite coarse decomposition of $M$ of multiplicity 1 (cyclic), although $M$ itself 
is not hyperfinite. For instance, one can take $\Gamma \curvearrowright B=R$ a non-commutative Bernoulli action, $B_0=B, B_1=\Bbb C$. 

We have thus shown: 

\proclaim{2.4. Corollary} There exists a non-Gamma $\text{\rm II}_1$ factor $M$ that 
contains hyperfinite subfactors $Q, P \subset M$ such that $Q\perp P$ and $\text{\rm sp}QP$ is $\|  \ \|_2$-dense in $M$.  
In other words, $M$ has an $R$-coarse pair decomposition of multiplicity $1$. 

\endproclaim

The next three results provide local criteria for a hyperfinite II$_1$ subfactor $Q\subset M$ (resp. for pairs of such factors $Q, P\subset M$) 
to be so that $_QL^2(M\ominus Q)_Q$ (resp. $_QL^2M_Q$) is factorial, or respectively properly infinite. 
They allow constructing $Q$ (resp. $Q, P$) recursively as inductive limits of finite dimensional factors. Note that the local conditions 2.5.1$^\circ$, $2.6.1^\circ$, 2.7.1$^\circ$ are ``Powers-type  
conditions'', in the spirit of [Po67]. 

For the purpose of the next statements, if $\xi$ is a unit vector in a Hilbert $\Cal H$ then we denote by $\omega_\xi$ the vector state 
on $\Cal B(\Cal H)$ implemented by $\xi$, i.e., $\omega_\xi(T)=\langle T(\xi), \xi \rangle$, $\forall T\in \Cal B(\Cal H)$. Also, 
if $\Cal H_0\subset \Cal H$ is a Hilbert subspace, then $p_{\Cal H_0}$ 
denotes the orthogonal projection of $\Cal H$ onto $\Cal H_0$. 

We denote $\text{\rm d}(p_{\Cal H_0}Tp_{\Cal H_0}, \Bbb Cp_{\Cal H_0})$ 
the distance in operator norm between the compression to  $\Cal H_0\subset \Cal H$ 
of an operator $T\in \Cal B(\Cal H)$ and the scalar multiples of $p_{\Cal H_0}$ (which is the identity in $\Cal B(\Cal H_0)$). 

It is useful to think of finite dimensional vector subspaces $\Cal H_0$ and the compressions $p_{\Cal H_0}Tp_{\Cal H_0}$ 
as ``windows through which we look at the operator $T$''. 

\proclaim{2.5. Proposition} Let $Q \subset M$ be a hyperfinite subfactor of the $\text{\rm II}_1$ factor $M$ and  
$\Cal L\subset L^2(M\ominus Q)$ a set of unit vectors that's dense in the set of unit vectors of $L^2(M\ominus Q)$. Then following conditions are equivalent: 

\vskip.05in

$1^\circ$ Given any finite set $F\subset \Cal L$, any finite 
dimensional subfactors $Q_0\subset Q$ and $\varepsilon >0$, 
there exists a finite dimensional subfactor $Q_1\subset Q$ that contains $Q_0$ 
such that $|\omega_\xi(XY)-\omega_\xi(X)\omega_\xi(Y)| \leq \varepsilon \|Y\|$ 
for all $Y\in (Q_1'\cap Q)\vee_{ \text{\rm Alg}}(Q_1'\cap Q)^{op}$, $X\in (Q_0\vee Q_{0}^{op})_1$, $\xi\in F.$ 

\vskip.03in

$2^\circ$ Given  any finite dimensional vector 
subspace $\Cal H_0\subset \text{\rm sp} \Cal L$ and $\varepsilon >0$ there exists a finite dimensional 
subfactor $ Q_1\subset Q$ such that if $Y\in ((Q_1'\cap Q)\vee_{ \text{\rm Alg}}(Q_1'\cap Q)^{op})_1$ 
then $\text{\rm d}(p_{\Cal H_0} Y p_{\Cal H_0}, \Bbb C p_{\Cal H_0})<\varepsilon$. 

\vskip.03in

$3^\circ$ Given  any finite dimensional vector 
subspace $\Cal H_0\subset L^2(M\ominus Q)$ and any $\varepsilon >0$, there exists a finite dimensional 
subfactor $ Q_1\subset Q$ such that if $Y\in ((Q_1'\cap Q)\vee (Q_1'\cap Q)^{op})_1$ 
then $\text{\rm d}(p_{\Cal H_0} Y p_{\Cal H_0}, \Bbb C p_{\Cal H_0})<\varepsilon$. 

\vskip.03in

$4^\circ$ $Q \vee Q^{op}\subset \Cal B(L^2(M\ominus Q))$ is a factor. 

\endproclaim
\noindent
{\it Proof}. $3^\circ \Leftrightarrow 4^\circ$. Writing $Q$ as the weak closure of an increasing sequence of finite dimensional 
subfactors $P_n \nearrow Q$, we have that $Q\vee Q^{op}\subset \Cal B(L^2(M\ominus Q))$ is a factor iff $(\cup_n P_n \vee P_n^{op})'\cap Q\vee Q^{op} = \Bbb C1$ 
in $\Cal B(L^2(M\ominus Q))$, 
a condition that's equivalent to the fact that $\forall F\subset L^2(M\ominus Q)$ finite, $\forall \varepsilon >0$, 
$\exists n$ such that if $T\in (\Cal B(L^2(M\ominus Q)))_1$ satisfies $[T, P_n]=0$, then $|\langle (T-\alpha_T)(\xi), \eta \rangle|< \varepsilon$, 
$\forall \xi, \eta \in F$, for some scalars $|\alpha_T|\leq 1$. Taking $F$ ``almost dense'' in a given finite dimensional vector 
space $\Cal H_0 \subset L^2(M\ominus Q)$, it follows that this last  condition is equivalent to $3^\circ$, by taking $Q_1=P_n$ for $n$ 
sufficiently large. 

$2^\circ \Leftrightarrow 3^\circ$ is immediate, by using the fact that $((Q_1'\cap Q)\vee_{Alg} (Q_1'\cap Q)^{op})_1$ is $so$-dense 
in $(Q_1'\cap Q)\vee (Q_1'\cap Q)^{op})_1$, Kaplansky density theorem, and the fact that $\overline{\text{\rm sp}\Cal L}=L^2(M\ominus Q)$. 

$1^\circ \Leftrightarrow 4^\circ$  We first make some observations. If $Q_0=P_0\subset P_1, ....\subset Q$ is an increasing  sequence of finite dimensional factors 
that exhaust $Q$ and denote $\tilde{P}_n=P_n \vee P_n^{op}$, $\tilde{P}=(\cup_n \tilde{P})_n)'' \subset \Cal B(L^2(M\ominus Q))$, then $4^\circ$ amounts to 
$\tilde{P}$ being a factor, i.e., $\tilde{P}'\cap \tilde{P}=(\cup_n \tilde{P}_n)'\cap \tilde{P}=\Bbb C$. Due to Kaplansky's density theorem,  $((Q_n'\cap Q)\vee_{Alg} (Q_n'\cap Q)^{op})_1$ is $so$-dense 
in $\tilde{P}_n'\cap \tilde{P}$, $\forall n$, so in particular any element in $(\Cal Z(\tilde{P}))_1$ is in the $so$-closure of $((Q_n'\cap Q)\vee_{Alg} (Q_n'\cap Q)^{op})_1$, $\forall n$. 

Since $\cap_n \tilde{P}_n'\cap \tilde{P}=\Cal Z(\tilde{P})$, it follows that $\Cal Z(\tilde{P})=\Bbb C$ iff any sequence $Y_n\in ((Q_n'\cap Q)\vee_{Alg} (Q_n'\cap Q)^{op})_1$ 
satisfies $\lim_{n\rightarrow\omega} Y_n\in \Cal C1$. Note that this is equivalent to the fact that $\forall \Cal H_0\subset L^2(M\ominus Q)$ finite dimensional 
(or just $\Cal H_0\subset \text{\rm sp}\Cal L$), $\forall \delta >0$, 
there exists $n$ large enough such that $\text{\rm d}(p_{\Cal H_0} Y p_{\Cal H_0}, \Bbb C p_{\Cal H_0})<\delta$, $\forall Y\in ((Q_n'\cap Q)\vee_{Alg} (Q_n'\cap Q)^{op})_1$. 

To see that $4^\circ \Rightarrow 1^\circ$, note that given $Q_0, F, \varepsilon$, there exists $\Cal H_0\subset L^2(M\ominus Q)$ finite dimensional 
that $\varepsilon/2$-contains $(Q_0\vee Q_0^{op})F(Q_0\vee Q_0^{op})$. If we assume $4^\circ$, then by the above observations there exists $n$ large enough such that any 
$Y\in \in ((P_n'\cap Q)\vee_{Alg} (P_n'\cap Q)^{op})_1$ there exists $\alpha_Y\in \Bbb C$ such that $\|p_{\Cal H_0} Y p_{\Cal H_0}- \alpha_Y p_{\Cal H_0}\| <\varepsilon/2$. 
Thus, for any such $Y$ and  any $X\in (Q_0\vee Q_0^{op})_1$, $\xi\in F$,  we have 
$$
|\omega_\xi(YX)-\omega_\xi(Y)\omega_\xi(X)|\approx_{\varepsilon} |\alpha_Y\omega_\xi(X)-\alpha_Y\omega_\xi(X)|=0.  
$$ 

Let us finally show $1^\circ \Rightarrow 4^\circ$.  From the above observations, we see that if $z\in (\Cal Z(\tilde{P})_+)_1$ and $F \subset L^2(M\ominus Q)$ finite, $\varepsilon >0$, are given, then 
there exists a finite dimensional factor $Q_0\subset Q$ large enough so that $(Q_0\vee Q_0^{op})_1$ contains a positive element $X$ that's $\varepsilon/4$-close to $z$ on 
$F$, as well as $n_0$ so that $\forall n\geq n_0$, 
$\exists Y\in ((Q_n'\cap Q)\vee_{Alg} (Q_n'\cap Q)^{op})_1$ that's $\varepsilon/4$-close to $z$ on $F$. Thus, $\omega_\xi(XY)\approx_{\varepsilon/2} \langle z^2\xi, \xi \rangle$ 
and $\omega_\xi(X)\omega_\xi(Y)\approx_{\varepsilon/2} (\langle z\xi, \xi \rangle)^2$, implying that $\langle z^2\xi, \xi \rangle\approx_{\varepsilon} (\langle z\xi, \xi \rangle)^2$. 
By taking $F$ sufficiently large and $\varepsilon >0$ small, this clearly implies 
$z$ must be a scalar. 
\hfill 
$\square$

\proclaim{2.6. Proposition} Let $Q, P \subset M$ be hyperfinite subfactors of the $\text{\rm II}_1$ factor $M$. 
Let also $\Cal L\subset L^2M$ be a set of unit vectors that's dense in the set of unit vectors of $L^2M$. Then following conditions are equivalent: 

\vskip.05in

$1^\circ$ Given any finite set $F\subset \Cal L$, any finite 
dimensional subfactors $Q_0 \subset Q, P_0\subset P$ and $\varepsilon >0$, 
there exist finite dimensional subfactors $Q_1\subset Q, P_1\subset P$ with $Q_1\supset Q_0$, $P_1\supset P_0$, 
such that $|\omega_\xi(XY)-\omega_\xi(X)\omega_\xi(Y)| \leq \varepsilon \|Y\|$ 
for all $Y\in (Q_1'\cap Q)\vee_{ \text{\rm Alg}}(P_1'\cap P)^{op}$, $X\in (Q_0\vee P_{0}^{op})_1$, $\xi\in F$. 

\vskip.03in

$2^\circ$ Given  any finite dimensional vector 
subspace $\Cal H_0\subset \text{\rm sp} \Cal L$, $\varepsilon >0$, there exist finite dimensional 
subfactors $ Q_1\subset Q$, $ P_1\subset P$, such that if $Y\in ((Q_1'\cap Q)\vee_{ \text{\rm Alg}}(P_1'\cap P)^{op})_1$ 
then $\text{\rm d}(p_{\Cal H_0} Y p_{\Cal H_0}, \Bbb C p_{\Cal H_0})<\varepsilon$. 

\vskip.03in

$3^\circ$ Given  any finite dimensional vector 
subspace $\Cal H_0\subset L^2(M)$ and any $\varepsilon >0$, there exist finite dimensional 
subfactors $ Q_1\subset Q$, $P_1\subset P$ such that if $Y\in ((Q_1'\cap Q)\vee (P_1'\cap P)^{op})_1$ 
then $\text{\rm d}(p_{\Cal H_0} Y p_{\Cal H_0}, \Bbb C p_{\Cal H_0})<\varepsilon$. 

\vskip.03in

$4^\circ$ $Q\vee P^{op}\subset \Cal B(L^2M)$ is a factor. 

\endproclaim
\noindent
{\it Proof}. The proof is the exactly the same as the proof of the equivalences in Proposition 2.5, so we leave the details as an exercise. 
\hfill 
$\square$ 

\proclaim{2.7. Proposition} Let $M$ be a separable $\text{\rm II}_1$ factor,  $Q, P\subset M$ subfactors and  
$\Cal L\subset L^2(M\ominus Q)$ $($respectively $\Cal L \subset L^2M)$ a  set of unit vectors that's dense 
in the set of unit vectors in $L^2(M\ominus Q)$ $($resp. $L^2M)$. The following conditions are equivalent: 

\vskip.05in

$1^\circ$ Given any finite set $F\subset \Cal L$, any finite 
dimensional subfactor $Q_0\subset M$ $($respectively finite dimensional  subfactors $Q_0, P_0\subset M)$ and $\varepsilon >0$, 
there exists a finite dimensional subfactor $Q_1\subset Q$ with $Q_1\supset Q_0$ $($respectively  finite dimensional subfactors $Q_1, P_1\subset M$, with $Q_1\supset Q_0$, 
$P_1\supset P_0)$ such that if one denotes by $a_\xi \in (Q_0'\cap Q_1) \vee (Q_0'\cap Q_1)^{op}$ $($resp. $a_\xi \in (Q_0'\cap Q_1)\vee (P_0'\cap P_1)^{op})$ the Radon-Nykodim 
derivative of $\omega_\xi$ with respect to $\tilde{\tau}=\tau \otimes \tau$ on $(Q_0'\cap Q_1)\vee (Q_0'\cap Q_1)^{op}$ $($resp. on $(Q_0'\cap Q_1)\vee (P_0'\cap P_1)^{op})$, 
then $\tilde{\tau}(e_{[\varepsilon, \infty)}(a_\xi)) \leq \varepsilon$.

\vskip.03in

$2^\circ$ $Q \vee Q^{op}\subset \Cal B(L^2(M\ominus Q))$ $($resp. $Q \vee P^{op}\subset \Cal B(L^2M))$ is properly infinite. 

\endproclaim
\noindent
{\it Proof}. We will only prove the case of a single hyperfinite factor $Q\subset M$ and leave the case of a pair of factors $Q, P\subset M$ as an exercise. 

Since $Q, Q^{op}$ are finite factors, the fact that $Q\vee Q^{op}$ has a finite part on $L^2(M\ominus Q)$ amounts to  having a unit vector $\xi \in L^2(M\ominus Q)$ 
such that on $p=[Q\xi Q^{op}]=[Q\vee Q^{op}(\xi)]\in (Q\vee Q^{op})'$ we have $Q\vee Q^{op}p\simeq Q\overline{\otimes} Q^{op}$. This in turn is equivalent 
to the fact that ${\omega_\xi}_{|Q\vee Q^{op}}=\tilde{\tau}(\cdot \ b)$ for some $b\in L^1(Q\overline{\otimes}Q^{op})_+$ with $\tilde{\tau}(b)=1$. 

If we write as before $Q$ as an inductive limit of some increasing sequence of finite dimensional factors $Q_0=P_0\subset P_1 \subset .... \nearrow Q$, 
this means that if we denote $a_n$ the expectation of $b$ onto $(P_0'\cap P_n)\vee (P_0'\cap P_n)^{op}$, then $\|a_n - a\|_{1, \tilde{\tau}}\rightarrow 0$, where 
$a=E_{(Q_0'\cap Q)\overline{\otimes} (Q_0'\cap Q)^{op}}(b)$. 

This also shows that having no such a finite unit vector is equivalent to the fact that, given any $\xi \in L^2(M\ominus Q)$, the Radon-Nykodim derivative $a_n$ of the restriction of the vector state $\omega_\xi$ 
to $(P_0'\cap P_n)\vee (P_0'\cap P_n)^{op}$ ``looses all its mass'', as $n\rightarrow \infty$, in other words $\tilde{\tau}(e_{[\varepsilon, \infty)}(a_n))\rightarrow 0$, $\forall \varepsilon >0$. 
This is also clearly equivalent to checking this property only for unit vectors $\xi$ in a dense subset $\Cal L$, thus proving $1^\circ \Leftrightarrow 2^\circ$. 

\hfill 
$\square$

\vskip.05in
\noindent
{\bf 2.8. Remarks}. $1^\circ$ One can show that the hyperfinite II$_1$ factor $R$ admits factorial  bimodule decompositions of any type I,  II$_1$, II$_\infty$, 
III$_\lambda$, $0\leq \lambda \leq 1$, both over a single subfactor, and over a pair 
of subfactors. Moreover, for coarse decompositions of $M=R$,  we can obtain any multiplicity in the 
``Jones' range'' $\{4\cos^2 \pi/n \mid n \geq 3\} \cup [4, \infty]$. However, 
as we mentioned before, it remains as an open problem whether the values left out, i.e., $(0, 1)\cup (1,4)\setminus \{4\cos^2 \pi/n \mid n\geq 3\}$, 
can be realized or not. The problem of finding the set of all possible values of multiplicitities of $R$-coarse decompositions of II$_1$ factors  
is very interesting, both for a specific factor $M$ (so viewed as an invariant and as symmetry question for $M$), 
and for the totality of all II$_1$ factors.  One would expect that 
any number that can be realized as the multiplicity of some $R$-coarse decomposition of a II$_1$ factor $M$, can also be realized 
as the multiplicity of coarse bimodule decomposition of $M\simeq R$. 

$2^\circ$ It may be true that if $M$ has a properly infinite $R$-bimodule decomposition, then it has properly infinite 
$R$-bimodule decompositions of any type. 

$3^\circ$  Given a II$_1$ factor $M$, the set of possible types that can appear as homogeneous (or even factorial) $R$-bimodule decompositions 
of $M$ (over a single copy of $R$, or over an $R$-pair) may be an interesting invariant to study for $M$. 

The only case where one knows to completely calculate this invariant are the free group factors and their amplifications, $M=L(\Bbb F_t)$, $1 < t \leq \infty$ 
(Dykema-Radulescu  interpolated free group factors [R92], [Dy93]), for which any homogeneous  $R$-bimodule decomposition is coarse  with infinite multiplicity.  More generally, 
one has the following consequence of 2.2 above and (Theorem 4.2 in [GP98]): 

\proclaim{2.9. Theorem} Let $M$ be one of the interpolated free group factors $L(\Bbb F_t)$, $1< t \leq \infty$. If $R_0, R_1\subset M$ are hyperfinite subfactors, 
then $_{R_0}L^2M_{R_1}$ has a coarse direct summand with infinite multiplicity. Equivalently, $R_0 \vee R_1^{op}$ has a finite direct summand $\simeq R_0\overline{\otimes} R_1^{op}$ 
with $(R_0\vee R_1^{op})'$ properly infinite. In particular, any homogeneous $R$-bimodule decomposition of $M$ is coarse with infinite multiplicity. 
\endproclaim
\noindent
{\it Proof}. If we assume $c=c(M; R_0, R_1) < \infty$ is finite, then by Corollary 2.2 we have $c(M^s; R_0^s, R_1^s)=c$ for any $s>0$. 
Choose  $s$ so that $n=(t-1)s^{-2}+1$ is an even integer satisfying $n/2-1 \geq c+2$. It follows that there exists a finite set $X\subset L^2M$ 
such that $|X|  \leq n/2 -1$ and $[R_0XR_1]=L^2M$. But this contradicts (Theorem 4.2 in [GP98]). 

\hfill 
$\square$

\heading 3.  Tight decomposition of II$_1$ factors 
\endheading

\noindent
{\bf 3.1. Definitions}. Let $M$ be a II$_1$ factor. A pair of subfactors $Q, P\subset M$ is {\it tight} if $Q\vee P^{op}=\Cal B(L^2M)$. If $P\subset M$ is a given 
subfactor and $Q\simeq R$ is a hyperfinite subfactor of  $M$ with the property that $Q, P\subset M$ is tight, then we also say that $Q$ is a $R$-tight complement of $P$. 

If $M$  has a tight decomposition $Q, P\subset M$ with both $Q$ and $P$ hyperfinite II$_1$  factors, then 
we say that $Q, P$ is an $R$-tight decomposition of $M$, and that $M$ is an  {\it $R$-tight  factor}  
(or simply a {\it tight} factor).

\proclaim{3.2. Proposition} Let $M$ be a $\text{\rm II}_1$ factor and $Q, P\subset M$ a pair of subfactors with the property that $\text{\rm sp}QP$ is $\| \ \|_2$-dense in $M$ and 
$(Q\cap P)'\cap M = \Bbb C1$. Then $(Q, P)$ gives a tight decomposition of $M$.
\endproclaim
\noindent
{\it Proof}. The argument is the same as the proof of (Proposition 2.2  in [GP98]), where this statement is being proved in the case $Q, P$ are hyperfinite. 
\hfill $\square$

\vskip.1in 
\noindent
{\bf 3.3. Examples}  $1^\circ$ If in Example 2.2.3$^\circ$  we take $B_0=B_1=B$ and $H\cap \Gamma_0$ to act freely on $B$ and ergodically on $\Cal Z(B)$, then 
$Q=B \rtimes H$, $P=B \rtimes \Gamma_0$  is a tight pair in $M$. Indeed, this is because these conditions imply that $Q\cap P=B \rtimes (H\cap \Gamma_0)$ has 
trivial relative commutant in $M=B \rtimes \Gamma$, while by Example 2.2.4$^\circ$ we already know that $QP$ is total in $M$, so Proposition 3.2 applies. 
If in addition $H$  is amenable and $B$ is AFD, then $Q$ is an $R$-tight complement of $P$ in $M$. So if $\Gamma_0$ is amenable as well,  
then $Q, P$ gives an $R$-tight decomposition of $M$.  

\vskip.05in

$2^\circ$ Let $Q\subset P$ be an extremal inclusion of $\text{\rm II}_1$ factors with finite index and $T \subset S$ its associated 
symmetric enveloping inclusion of $\text{\rm II}_1$ factors, as defined in [P94]. Thus, $S$ is generated by commuting copies 
of $P, P^{op}$ and by a projection $e$ of trace $[P:Q]^{-1}$ that implements at the same time the expectation of $P$ onto $Q$ 
and of $P^{op}$ onto $Q^{op}$, while $T$ is generated by $P, P^{op}$. Then there exists a choice of 
a tunnel-tower for $Q\subset P, Q^{op}\subset P^{op}$,  as in [P97], such that the corresponding enveloping factors $P_\infty, P^{op}_\infty\subset S$ 
have the property that $(P_\infty \cap P^{op}_\infty)'\cap  S = \Bbb C$. By Proposition 3.2, this shows that $P_\infty, P^{op}_\infty$ provide a tight decomposition of  $S$, which is 
$R$-tight when $P$ is hyperfinite. 

Also, if $P$ is hyperfinite and the standard graph $\Gamma_{Q\subset P}$ is ergodic (i.e., $P'\cap P_\infty$ is a factor, for instance if $\Gamma_{Q\subset P}=A_\infty$), 
then $P_\infty$ and $T=P\vee P^{op}$ provide 
a tight decomposition of $S$ as well, which is hyperfinite whenever $P\simeq R$. 

\vskip.05in

$3^\circ$ More generally, let $Q\subset P$ and $\tilde{Q} \subset \tilde{P}$ be extremal inclusions of factors with the same standard invariant, 
$\Cal G_{Q\subset P}=\Cal G_{\tilde{Q}\subset \tilde{P}}$, and denote by $T\subset S$ its associated enveloping inclusion of II$_1$ factors, 
obtained by taking the ``concatenated product'' of the two subfactors, 
as defined in (Remark 2.5.1$^\circ$ in [P97]). Thus, $S$ is generated by commuting copies of $\tilde{P}, P^{op}$ 
and a projection $e$ of trace $\lambda=[P:Q]^{-1}=[\tilde{P}: \tilde{Q}]^{-1}$ implementing both the expectation of $\tilde{P}$ onto $\tilde{Q}$ 
and of $P^{op}$ onto $Q^{op}$, and with $T=\tilde{P}\vee P$. Then there exists a choice of Jones $\lambda$-projections $\{e_n\}_{n\in \Bbb Z}\subset T$ 
with $e_0=e$, $e_n\in \tilde{P}$ for $n<0$ and $e_n \in P^{op}$ for $n>0$, such that: $\tilde{P}_n=\tilde{P}\cap \{e_0, ..., e_n\}'$, $n=-1, -2, ...$, 
gives a tunnel for $\tilde{Q}\subset \tilde{P}$ with $\tilde{P}_{-1}=\tilde{Q}$; $P^{op}_{-n}=P^{op}\cap \{e_0, ..., e_n\}'$, $n=1, 2, ....$, 
gives a tunnel for $Q^{op}\subset P^{op}$ with $P^{op}_{-1}=Q^{op}$. Like in the case of usual symmetric enveloping 
inclusion of a subfactor, there exist choices of $\{e_n\}_n$ 
such that  $P_\infty, P^{op}_\infty$ is a tight  pair in $S$ (due to the fact that $P_\infty P^{op}_\infty$ is total in $S$ and 
$(\tilde{P}_\infty \cap P^{op}_\infty)'\cap S=\Bbb C$, and using Proposition 3.2).  
If in addition $\Gamma_{Q\subset P}=\Gamma_{\tilde{Q}\subset \tilde{P}}$ is ergodic and $P$ is hyperfinite, then $P^{op}_\infty$ is 
an $R$-tight complement of $T$ in  $S$.

\proclaim{3.4. Proposition}   Let $M$ be a separable $\text{\rm II}_1$ factor and $\Cal L \subset L^2M$ a countable set of 
unit vectors that's dense in the set of all unit vectors of $L^2M$. The following conditions are equivalent: 
\vskip.05in 

$(a)$ Given any finite dimensional subspace $\Cal H_0\subset \text{\rm sp} \Cal L$, any finite dimensional subfactors $Q_0, P_0\subset M$ and $\varepsilon >0$, there exist finite dimensional 
subfactors $Q_1, P_1 \subset M$ such that $Q_1\supset Q_0$, $P_1\supset P_0$ and such that if $Y\in ((Q_1'\cap M)\vee_{ \text{\rm Alg}}(P_1'\cap M)^{op})_1$ 
then $\text{\rm d}(p_{\Cal H_0} Y p_{\Cal H_0}, \Bbb C p_{\Cal H_0}) < \varepsilon$. 
\vskip.03in

$(b)$ Given any finite subset $F$ of vector  states implemented by vectors in $\Cal L$, any finite dimensional subfactors $Q_0, P_0\subset M$ and $\varepsilon >0$, 
there exist finite dimensional 
subfactors $Q_0\subset Q_1 \subset M$, $P_0\subset P_1\subset M$ such that $\|E_{(Q_1\vee P^{op}_1)'}(\varphi - \psi)\|\leq \varepsilon$, for all $\varphi, \psi\in F$.  

\vskip.03in 

$(c)$ For any finite dimensional factors $Q_0, P_0 \subset M$, any $\eta\in \Cal L$ and 
any $\varepsilon>0$, there exist finite dimensional factors $Q_0 \subset Q_1 \subset M$, $P_0 \subset P_1\subset M$ and elements 
$x_j\in Q_1, y_j\in P_1$ such that $\|\Sigma_j x_j \otimes_{min} y_j^{op}\|\leq 1$ and $\|\Sigma_j x_j y_j -\eta\|_2\leq \varepsilon$. 

\vskip.03in
$(d)$ $M$ is $R$-tight.  
\endproclaim
\noindent
{\it Proof}. Let $\Cal H_n \subset \text{\rm sp} \ \Cal L$ be an increasing sequence of finite dimensional subspaces that exhaust $\text{\rm sp} \ \Cal L$ 
(thus $\cup_n \Cal H_n$ is dense in $L^2M$). Condition $(a)$ allows constructing recursively an increasing sequence of finite dimensional factors $Q_n, P_n \subset M$ 
such that  any $T\in \Cal B(L^2M)$ that commutes with all $Q_n, P_n^{op}$ must be arbitrarily close  to a scalar multiple of $1$ on 
$\Cal H_n \nearrow L^2M$. Thus, if $Q, P$ denote their corresponding limit, then $(Q\vee P^{op})'\cap \Cal B(L^2M)=\Bbb C 1$, 
showing that $(a)\Rightarrow (d)$. 

Similarly,  since the set $S(\Cal L)$ of vector states implemented by $\Cal L$ is a countable set of states on $\Cal B(L^2M)$ that's dense in the space of vector states, condition  
$(b)$ can be used to construct recursively an increasing sequence of finite dimensional factors $Q_n, P_n \subset M$ 
such that $\lim_n \|E_{(Q_n \vee P^{op}_n)'}(\varphi - \psi)\|=0$, $\forall \varphi, \psi \in S(\Cal L)$. By (Corollary 2.5 in [P19]), 
this proves $(b)\Rightarrow (d)$. 

To see that $(d)$ implies both $(a)$, $(b)$, let $Q, P\subset M$ be hyperfinite subfactors satisfying $(Q\vee P^{op})'=\Bbb C1$. 
Since this latter condition is invariant to conjugating $Q, P$ by unitary elements, we may assume $Q, P$ contain some prescribed 
finite dimensional factors $Q_0\subset Q, P_0 \subset P$. If $Q_0 \subset Q^n \nearrow Q$, $P_0 \subset P^n \nearrow R$ are 
increasing sequences of finite dimensional subfactors generating $Q, P$, then taking $Q_1$, respectively  $P_1$ to be $Q^n$, 
respectively $P^n$, for $n$ sufficiently large, shows that $(d) \Rightarrow (a)$, $(d) \Rightarrow (b)$. 

The fact that $(d)$ implies $(c)$ is trivial by Kaplansky's density theorem. Conversely, if $(c)$ is satisfied, and one choses a 
sequence $\Cal L$ of unit vectors in $L^2M$ that's dense in the set of unit vectors of $L^2M$,  then one can construct recursively 
a pair of increasing sequences of finite dimensional factors $Q^n$, respectively $P^n$ in $M$ such that if one denotes $Q, P$ their respective limits, 
then for any $\xi \in \Cal L$ and $\varepsilon >0$ there exists $T\in (Q\vee P^{op})_1$ such that $\|T(\hat{1})-\xi\|\leq \varepsilon$. If $\eta \in L^2M$ is an arbitrary unit vector 
that's a limit of some $\xi_n \in \Cal L$, then let $T_n \in (Q\vee P^{op})_1$ be so that $\|T_n(\hat{1})-\xi_n\|\leq 2^{-n}$. Thus, if $T$ is a weak limit of $T_n$, 
then $T\in (Q\vee P^{op})_1$ will satisfy $\langle T(\hat{1}), \eta \rangle =1$. Together with the fact that $\|T\|\leq 1$ and $\|\eta\|=1$, this implies $T(\hat{1})=\eta$. 
It also implies $T^*(\eta)=\hat{1}$. 

Thus, $Q\vee P^{op}(\eta)=L^2M$ for any $\eta\neq 0$, showing that $Q\vee P^{op}=\Cal B(L^2M)$. 
\hfill $\square$

\heading 4.  Stably single generated (SSG) II$_1$ factors  \endheading

In this section we make some general remarks about the (minimal) number of generators of a II$_1$ factor and 
about the stably single generated property. All of these results can be readily deduced from ([Sh05], [DSSW07]), 
but we will provide direct arguments, for the reader's convenience. 

We begin by noticing a scaling formula for the number of generators, which is more or less implicit in ([V88]): 

\proclaim{4.1. Lemma} $1^\circ$ If a $\text{\rm II}_1$ factor $M$ can be generated by $n\geq 2$ selfadjoint elements and $k\geq 1$ is an integer, then 
$M^{1/k}$ can be generated by $(n-1)k^2+1$ selfadjoint elements. 

$2^\circ$  If a $\text{\rm II}_1$ factor $P$ can be generated by 
$(n-1)k^2+1$ selfadjoint elements, for some integers $n\geq 2, k\geq 1$, then $\Bbb M_{k}(P)$ can be generated by $n$ selfadjoint elements. 

$3^\circ$ If $M$ is a $\text{\rm II}_1$ factor and $p\in M$ a non-zero projection such that $pMp$ is generated by $n$ self-adjoint elements, 
then $M$ can be generated by $n$ self-adjoint elements.  

\endproclaim
\noindent
{\it Proof}. Let us first notice that if $a_1, ..., a_n \in M$ are self adjoint elements generating a II$_1$ factor $M$, then we can replace 
 them by positive invertible elements $b_1, ..., b_n$ so that $\{b_i\}''$ is a MASA, $\forall i$.

$1^\circ$. Let $a_1, ..., a_n\in M_h$ be positive invertible 
elements generating $M$, for some $n\geq 2$. Let $e_{ii}, 1\leq i \leq k$, be a partition of $1$ with projections 
of trace $1/k$ that commute with $a_1$. 

For each $2\leq j \leq k$, by using polar decomposition  we can write the element $e_{11}a_2e_{jj}$, 
in the form  $a^2_je_{1j}$ where $e_{1j}$ is a partial isometry with left support $e_{11}$ and right support $e_{jj}$, and $a^2_j$ is a positive element 
in $e_{11}Me_{11}$. 

Denote $e_{j1}=e_{1j}^*$. Since each single element $e_{1i}a_me_{j1}$ 
with $k\geq j > i\geq 2$, $n\geq m\geq 3$, is the sum between its real part and $i$-times its imaginary part, counting all the resulting selfadjoint elements in $e_{11}Me_{11}$ we get: 

$(a)$ $kn$ elements of the form $e_{1j}a_me_{j1}$, $1\leq m \leq n$, $1\leq j \leq k$; 

$(b)$ $k-1$ elements of the form $a^2_j$, $2\leq j \leq k$; 

$(c)$ $2$ elements for each $e_{1i}a_2e_{j1}$ $k\geq j>i \geq 2$, for a total of $2(k-1)(k-2)/2=(k-1)(k-2)$ elements;  

$(d)$ $2$ elements for each $e_{1i}a_me_{j1}$ $k\geq j>i \geq 1$, $n\geq m \geq 3$, for a total of 
$2(n-2) k(k-1)/2=(n-2)k(k-1)$ elements. Summing up, we get an overall total of 
$$
nk+(k-1)+(k-1)(k-2)+ (n-2)k(k-1)
$$
$$
=nk +(k-1)^2 + nk^2-nk-2k^2+2k
$$
$$
=nk^2 -k^2 +1=(n-1)k^2+1  
$$
elements. So the statement follows once we notice that any element in $e_{11}Me_{11}$ is a linear combination of products of the form $e_{11}a_{i_1}a_{i_2}.... a_{i_t}e_{11}$ and that we can 
write each product  $a_sa_{s'}$ in the form $1a_{s}1a_{s'}1=\sum_{i,j,k} e_{ii}a_{s}e_{jj}a_{s'}e_{kk}$, with $e_{ii}=e_{i1}e_{i1}^*$, making each element in $e_{11}Me_{11}$ 
be a linear combination of elements of one of the forms $(a)-(d)$. This shows that $e_{11}Me_{11}$ is generated by $(n-1)k^2+1$ elements.

$2^\circ$ The proof of this part is very similar to $1^\circ$. We detail it in the case $n=2$ and leave the general case as an exercise. 

So let $\{e_{ij} \mid 1\leq i,j \leq k\}$ be matrix units for $\Bbb M_k(\Bbb C) \subset \Bbb M_k(P)$. We label the $k^2+1$ generators 
of $P$ as $\{a_i \mid 1\leq i \leq k\}$, $\{b^i_j \mid 1\leq i \leq j \leq k\}$, $\{u^{s}_t \mid 2\leq s < t \leq k\}$, where the first two sets are made of positive invertible 
elements while the last one is made of unitary elements (which can be viewed as exponentials of selfadjoint elements). Thus, the total number is indeed 
$$
k+k(k+1)/2 + (k-1)(k-2)/2$$
$$
=k+(k^2+k+k^2-3k+2)/2=k^2+1.
$$ 

Now let $a\in  M:=\Bbb M_k(P)$ be a positive invertible element generating a MASA that contains $\{a_ie_{ii}, e_{ii} \mid 1\leq i \leq k\}$ and 
$$
b=\sum_{1\leq i \leq k} b^i_ie_{ii} + \sum_{1\leq i < j \leq k} (b^i_ju^i_j e_{ij} + e_{ji}{u^i_j}^*b^i_j). 
$$

Then clearly $a, b$ are selfadjoint elements in $M$ and if one denotes by $M_0$ the von Neumann subalgebra of $M$ generated by $a, b$ then 
 $e_{1j} \in M_0$, $1\leq j \leq k$, and $e_{11}M_0e_{11}=Pe_{11}$. 

Indeed, we have  $e_{ii}\in \{a\}''\subset M_0$,  
so $e_{11}be_{jj}\in M_0$ as well, implying that $e_{1j} = (e_{11}be_{jj}b^*e_{11})^{-1/2}(e_{11}be_{jj})\in M_0$. This also shows that $a_ie_{11}= 
e_{1j}ae_{j1}, b_i=e_{1j}be_{j1}\in e_{11}M_0e_{11}$, $\forall i$. Also, for each $b^i_ju^i_j e_{11}=e_{1i}be_{j1}\in e_{11}M_0e_{11}$. 
Thus, $b^i_j$ and $u^i_j$, which give the positive and unitary elements in the polar decomposition of $b^i_ju^i_j e_{11}$, 
belong to $e_{11}M_0e_{11}$, for all $1\leq i < j \leq n$. 

Thus, $Pe_{11}\subset e_{11}M_0e_{11} \subset e_{11}Me_{11}=Pe_{11}$, while $e_{1j}\in M_0$ as well, 
which altogether implies $M_0=M$. 

$3^\circ$ It is sufficient to prove this for $\tau(p)\geq 1/2$. Let $a_1, ..., a_n\in pMp$ be selfadjoint elements generating $pMp$ as a von Neumann algebra. 
We may clearly assume that each one of the weakly closed $^*$-algebras generated by $a_j$, $j=2, ..., n$, contains $p=1_{pMp}$. 
Let $v\in M$ be a partial isometry such that 
$v^*v=p$, $vv^*\leq 1-p$. Denote $b_1=a_1 + v+v^* \in M$. Then $b_1^*=b_1$ and since the weakly close $^*$-algebra generated by $a_2$ contains $p$, 
it follows that the weakly closed $^*$-algebra $N\subset M$ generated by $b_1, a_2, ..., a_n$ contains $pb_1p=a_1$, thus $pMp\subset N$. 
Since we also have $v=b_1p-pb_1p \in N$, it follows that $M=N$. 
\hfill $\square$

\vskip.1in
\noindent
{\bf 4.2. Notation} If $M$ is a von Neumann algebra, then we denote by $\text{\rm ng}(M)$ the minimal number $n \in \Bbb N \cup \{\infty\}$ with the property that 
$M$ can be generated by $n$ selfadjoint elements. Note that if $M$ is a II$_1$ factor, one actually has ng$(M)\geq 2$.

\proclaim{4.3. Corollary} Let  $M$ be a $\text{\rm II}_1$ factor. 

$1^\circ$ The function $\text{\rm ng}(M^t)$ is non-increasing in $t>0$. 

$2^\circ$ Let $k\geq 1$ be an integer. Then $\text{\rm ng}(M) \leq n$ if and only if $\text{\rm ng}(M^{1/k}) \leq (n-1)k^2+1$. 

\endproclaim
\noindent
{\it Proof}. Part $1^\circ$ is trivial from 4.1.3$^\circ$ while $2^\circ$ is an immediate consequence of $4.1.1^\circ, 4.1.2^\circ$.  
\hfill $\square$

\proclaim{4.4. Corollary} Let  $M$ be a $\text{\rm II}_1$ factor. 

$1^\circ$ The factors $M^t$ are either all finitely generated $($i.e., $\text{\rm ng}(M^t)<\infty$, $\forall  t>0)$, 
or all infinitely generated $($i.e., $\text{\rm ng}(M^t)=\infty$, $\forall t>0)$.  

$2^\circ$ If $M^t$ are finitely generated then there exists a constant $K=K(M)$ such that $ \text{\rm ng}(M^t)\leq K t^{-2}$, $\forall 0<t\leq 1$. 
\endproclaim
\noindent
{\it Proof}.  Both properties are immediate consequences of Corollary 4.3. 
\hfill $\square$

\vskip.05in
\noindent
{\bf 4.5. Definition} A II$_1$ factor $M$ is {\it stably single generated}, or has the {\it stable single generation}   ({\it SSG}) property, 
if $M^t$ is single generated (in other words, $\text{\rm ng}(M^t)=2$) for all $t>0$. 

\vskip.05in
\noindent
{\bf 4.6. Remark}. Voiculescu's methods for measuring the contribution/thinness of sets of generators of a II$_1$ factor $M$ in (Section 7 of [V96]), 
where used by J. Shen in ([Sh05]) to isolate a numerical invariant denoted  $\Cal G(M)$, which  provides a ``weigthed counting'' 
of the minimal number of generators  of $M$, minus 1, and lies in $[0, \infty]$.  This invariant was further studied in [DSSW07]. It follows from results in these two papers that the 
SSG property of $M$ in our sense  is equivalent to $\Cal G(M)=0$. Thus, results about the stability/permanence of the property $\Cal G(M)=0$ in ([Sh05], [DSSW07]) 
provide stability/permanence properties for SSG of $M$. We enumerate below just a few such properties, but which we deduce directly from Corollary 4.3 and  ([GP98]). 
We refer the reader to (5.16  in [Sh05] and [DSSW07]) for more permanence results and examples of classes of factors  $M$ with $\Cal G(M)=0$, which are thus SSG. 

\proclaim{4.7. Corollary} Let $M$ be a $\text{\rm II}_1$ factor. If $\liminf_{t\rightarrow 0} \ t^2\text{\rm ng}(M^t)=0$, then  $M$ is SSG. 
In particular, if the non-increasing function $t \mapsto \text{\rm ng}(M^t)$ is uniformly bounded, then $M$ is SSG. 
\endproclaim
\noindent
{\it Proof}. This is trivial by Corollary 4.3. 
\hfill $\square$

\proclaim{4.8. Corollary} Let $M$ be a $\text{\rm II}_1$ factor. If $M$  is finitely generated and has non-trivial fundamental group, 
then $M$ is SSG.  
\endproclaim
\noindent
{\it Proof}. This is clear by Corollary 4.7. 

\hfill $\square$

For the next result, recall from [P16] that a MASA $A$ in a II$_1$ factor is an {\it s-MASA} if $A \vee A^{op}$ is a MASA in $\Cal B(L^2M)$, and that a factor is {\it s-thin} 
if it has an s-MASA.  Note that by [FM77], any Cartan MASA is an s-MASA. 

\proclaim{4.9. Proposition}  Any separable $\text{\rm II}_1$ factor $M$ satisfying one of the properties below is  SSG: 

$1^\circ$ $M$ is non-prime, i.e.,  $M\simeq N\overline{\otimes} P$ with both $N, P$ of type $\text{\rm II}_1$;  

$2^\circ$ $M$ has  a Cartan subalgebra, more generally an s-MASA; 

$3^\circ$  $M$ has the property Gamma. 

$4^\circ$ $M$ is $R$-tight, more generally $M$ has a properly infinite $R$-bimodule decomposition. 

$5^\circ$ $M$ is weakly $R$-thin, i.e., there exist a pair of hyperfinite factors $R_0, R_1\subset M$ and a finite set $X\subset L^2M$ 
such that $[R_0 X R_1]=L^2M$. 
\endproclaim 
\noindent
{\it Proof}. The fact that $M$ is either non-prime, or has the property Gamma, or is tight, or has a Cartan subalgebra, are all properties that are obviously stable to 
amplifications by arbitrary $t>0$. So $1^\circ$, $3^\circ$, $4^\circ$ follow from the fact that a factor having any of these properties is single generated by (6.2 in [GP98]). 
The property of being s-thin (i.e., having an s-MASA) has been shown to be stable to amplifications by arbitrary $t>0$ in (3.8 of [P16]), so $2^\circ$ follows from (3.5 in [P82]). 

To see that $5^\circ$ holds true, let $c(M; R_0, R_1)$ denote the coupling constant of $(R_0\vee R_1^{op})_{fin}$ which is finite because it is bounded by $|X|$. 
Note that ng$(M)\leq \text{\rm ng}(R_0)+\text{\rm ng}(R_1) + 2 (|X| +1) = 2c(M; R_0, R_1) + 6 < \infty$. Similarly, 
by considering the $R_0^t-R_1^t$ bimodule decomposition of $M^t$, one gets ng$(M^t)\leq  2c(M^t; R_0^t, R_1^t) + 6$. But by Lemma 2.1, 
we have $c(M; R_0, R_1)=c(M^t; R_0^t, R_1^t)$. This shows that ng$(M^t)$ is uniformly bounded, so by Corollary 4.7, $M$ is SSG. 
\hfill $\square$

\proclaim{4.10. Proposition}  $(a)$ If $N$ is an SSG $\text{\rm II}_1$ factor and $\Gamma \curvearrowright^\sigma N$ 
is a free cocycle action of a countable group, then $M=N \rtimes \Gamma$ is SSG. 

$(b)$ If $N\subset M$ is an irreducible, quasiregular inclusion of separable $\text{\rm II}_1$ factors and $N$ is SSG, then $M$ is SSG. 

$(c)$ If $N \subset M$ is an inclusion of $\text{\rm II}_1$ factors with finite index then $M$ is SSG iff $N$ is SSG. 

\endproclaim 
{\it Proof}. Part $(a)$ amounts to showing that if $N\subset M$ is an irreducible regular inclusion of separable II$_1$ factors and $N$ is SSG, then $M$ is SSG. 
But if $N\subset M$ satisfies these hypothesis, then so does $(N\subset M)^t = (N^t \subset M^t)$, for any $t>0$. 

Thus, by Corollary 4.7, in order to prove $(a)$ it is sufficient to show that 
under the given conditions $M$ follows  generated by two  elements (i.e., $\text{\rm ng}(M)\leq 4$). So let 
$\{u_g \mid g\in \Gamma\}\subset M$ be the canonical unitaries implementing the $\Gamma$-action $\sigma$ and denote $b=\sum_g 2^{-n_g} u_g$, where $\{n_g\}_g$ are distinct positive integers. 

We claim that $M$ is generated by $N$ and $b$, thus by two elements. Indeed, one has that the von Neumann algebra $M_0$ generated by 
$N$ and $b$ contains the element of minimal $\| \ \|_2$-norm in the weak closure of the convex hull of $\{ub\sigma_g(u)^* \mid u\in \Cal U(N)\}$, which is clearly equal 
to $2^{-n_g}u_g$. Thus, $M_0$ contains $N$ and $\{u_g\}_g$, so $M_0=M$. 

Part $(b)$, which generalizes $(a)$, is proved in exactly the same way, using the fact that $N\subset M$  irreducible and quasi-regular is a property that's stable 
to amplification by any $t>0$, and then by a fact from [P01], showing that if $N\subset M$ is irreducible and quasi-regular then there exist partial isometries 
$v_n \in M$ such that as $N-N$ Hilbert bimodules, we have $L^2M = \oplus L^2(N v_n N)$, with $v_n^*v_n \in N$ and $v_nNv_n^*=Q_n q_n$, where $Q_n\subset N$ is a finite index subfactor 
with $q_n\in Q_n'\cap M$. 

To prove $(c)$, note that by [PiP84], $M$ is generated by $N$ and by an orthonormal basis $M$ over $N$ that has less than $[M:N]+1$ elements. Since the index is stable to amplification, 
it follows that $M^t$ can be generated by less than $[M:N] +2$ elements, $\forall t>0$, so $N$ SSG implies $M$ is SSG, by Corollary 4.7. This also shows that if $M$ is SSG 
then $\langle M, e_N \rangle \simeq N^{[M:N]}$ is SSG, so $N$ follows SSG as well. 
\hfill $\square$

\heading 5.  Local characterization of SSG 
\endheading

Let us first notice a general ``weak spectral gap'' type characterization of the fact that a II$_1$ factor $M$ is generated by a 
finite set of unitaries $\{u_i\}_i$. 

\proclaim{5.1. Lemma} Let $M$ be a von Neumann algebra represented on the Hilbert $\Cal H$ and  $u_1, ..., u_n \subset M$ 
a finite set of unitary elements. Let $\Cal M\subset \Cal B(\Cal H)$ be a von Neumann algebra that contains $M$. 
Consider the conditions:  

\vskip.05in 
$1^\circ$ $\{u_i\}_i$  generate $M$ as a von Neumann algebra. 

$2^\circ$ Given any $\varepsilon > 0$ and any 
finite $F\subset \Cal H$, there exist $\delta>0$ and a finite $E\subset \Cal H$ such that if $T\in (\Cal M)_1$,  
satisfies $|\langle (Tu_i-u_iT)(\eta), \eta' \rangle | \leq \delta$, $\forall \eta, \eta'\in E$, then there exists $T'\in (M'\cap \Cal M)_1$ 
with  $|\langle (T-T')\xi, \xi' \rangle | \leq \varepsilon$, $\forall \xi, \xi' \in F$.  

\vskip.05in 

Then $1^\circ \Rightarrow 2^\circ$. If in addition $\Cal M=\Cal B(\Cal H)$, then $2^\circ \Rightarrow 1^\circ$ as well. 
  
\endproclaim 
\noindent
{\it Proof}. $1^\circ \Rightarrow 2^\circ$ Let $J$ be the set of all finite subsets of $\Cal H$, which we view as a directed set with respect to the order 
given by inclusion. Assume by contradiction that $2^\circ$ fails to be true. 
This means there exists $\varepsilon_0>0$ and a finite set $F_0\subset \Cal H$, such that for any $j\in J$ there exists $T_j\in (\Cal M)_1$ 
satisfying  
$$
\max \{|\langle (T_j u_i-u_iT_j)\eta, \eta' \rangle| \mid \eta, \eta'\in j\}  \leq |j|^{-1}, \forall 1\leq i  \leq n, \tag 5.1.1 
$$
but  
$$
\Sigma_{\xi, \xi'\in F_0} |\langle (T_j -T')\xi, \xi' \rangle|  \geq \varepsilon_0,  \forall T'\in (M'\cap \Cal M)_1. \tag 5.1.2 
$$

Let $T\in (\Cal M)_1$ be a Banach (weak) limit of the net $(T_j)_j$. By $(5.1.1)$, it follows that $Tu_i=u_iT$, $\forall i$. Since $u_i$ 
generate $M$ as a von Neumann algebra, this implies $T\in(M'\cap \Cal M)_1$. But if we take  in $(5.1.2)$ 
the Banach limit $T$ of the net $(T_j)_j$,  
then we get 
$$
\Sigma_{\xi, \xi'\in F_0} |\langle (T -T')\xi, \xi' \rangle|  \geq \varepsilon_0,  \forall T'\in (M'\cap \Cal B(\Cal H))_1, \tag  5.1.3 
$$
a contradiction. 

In the case $\Cal M=\Cal B(\Cal H)$, by von Neumann's bicommutant theorem the fact that $\{u_i\}_i$ generate $M$ 
as a von Neumann algebra means that the 
commutant of $\{u_i\}_i \cup \{u_i^*\}_i$ in $\Cal M=\Cal B(\Cal H)$ coincides with 
$M'\cap \Cal B(\Cal H)$. It is clear that $2^\circ$ implies this latter fact, showing that $2^\circ \Rightarrow 1^\circ$. 

\hfill 
$\square$ 

\vskip.05in 

We next notice that SSG factors can be generated by a pair of triadic matrix algebras with the same diagonal. 
This will provide a good intuition for our conjectures in the next section. 

\proclaim{5.2. Lemma}  If $P$ is generated by two selfadjoint elements then $M=\Bbb M_3(P)$ is generated by two matrix subalgebras 
$B_0, B_1 \simeq \Bbb M_3(\Bbb C)$, that have a common diagonal. 
\endproclaim
\noindent
{\it Proof}. Assume $P$ is generated by the unitary elements $u_1, u_2$. Take $B_0=\Bbb M_3(\Bbb C)\subset \Bbb M_3(P)=M$ with a system of matrix units 
$\{e_{ij}\mid 0\leq i, j\leq 2\}$.  Then define $B_1$ to be the $^*$-algebra generated by $e_{00}, e_{11}, e_{22}$ and $u_1e_{01}, u_2e_{02}$. 
Then $B_1\simeq \Bbb M_3(\Bbb C)$  has the same diagonal as $B_0$, and clearly $B_0, B_1$ generate $M$. 
\hfill $\square$

\proclaim{5.3. Corollary}  Let $M$ be a $\text{\rm II}_1$ factor. The following conditions are equivalent: 

\vskip.05in 
$1^\circ$ $M$ is SSG; 

\vskip.05in

$2^\circ$ $M^t$ is generated by two matrix subalgebras 
$B_0, B_1 \simeq \Bbb M_3(\Bbb C)$, that have a common diagonal, for any $t>0$. 

\vskip.05in 

$3^\circ$ $M$ contains two sequences of matrix units $\{e^{0,n}_{ij}\}_{0\leq i,j \leq 2}$, $\{e^{1,n}_{ij}\}_{0\leq i,j \leq 2}$, $n\geq 1$, such that for each $n\geq 1$ 
we have: 
\vskip.03in
$(a)$ $e^{0,n}_{ii}=e^{1,n}_{ii}$, $\forall i$; 

$(b)$  $\sum_{i=0}^2 e^{0,n}_{ii}=e^{0,n-1}_{00}$, where 
$e^{0,0}_{00}=1$; 

$(c)$ $\{e^{0,n}_{ij}\}_{0\leq i,j \leq 2}$, $\{e^{1,n}_{ij}\}_{0\leq i,j \leq 2}$ generate $e^{0,n}_{00}Me^{0,n}_{00}=e^{1,n}_{00}Me^{1,n}_{00}$. 

\vskip.05in

$4^\circ$ $M$ contains two hyperfinite subfactors $R_0, R_1\subset M$ with the property that $R_0\cap R_1$ contains a decreasing sequence of projections $\{e_n\}_n$ 
with $\tau(e_n)=3^{-n}$ such that $e_nR_0e_n, e_nR_1e_n$ generate $e_nMe_n$ for all $n$. 
\endproclaim
\noindent
{\it Proof}. This is now immediate by Lemma 5.2. 
\hfill $\square$

\vskip.05in

The next result shows permanence of SSG to inductive limits. As we mentioned before, since SSG is equivalent to the invariant 
in ([Sh05]) satisfying $\Cal G(M)=0$, this result is just (5.16.$(iv)$ in [Sh05]). However, we provide a short direct proof which utilizes  the 
generation of $M$ by pairs of triadic factors, as in 5.2, 5.3 above, which we think is relevant for the discussions in Section 6. 

\proclaim{5.4. Theorem} If $N_n \nearrow M$ are $\text{\rm II}_1$ factors and $N_n$ is SSG, $\forall n$, then $M$ is SSG. 
\endproclaim 
\noindent
{\it Proof}. We denote $e^{0,0}_{00}=e^{1,0}_{00}=1$ and construct recursively two sequences of  $3 \times 3$ 
matrix units $\{e^{0,n}_{ij}\}_{0\leq i,j \leq 2} \subset N_n$, $\{e^{1,n}_{ij}\}_{0\leq i,j \leq 2}\subset N_n$, $n \geq 1$ with the following properties: 

\vskip.05in 

$(a)$ $e^{0,n}_{ii}=e^{1,n}_{ii}$  for all $0\leq i \leq 2$ and all $n\geq 1$.

\vskip.05in

$(b)$ $\sum_{i=0}^2 e^{0, n}_{ii}=e^{0,n-1}_{00}=e^{1, n-1}_{00}=\sum_{j=0}^2 e^{1, n}_{jj}$, for all $n\geq 1$;

\vskip.05in 

$(c)$  $\{e^{0,n}_{ij}\}_{0\leq i,j \leq 2} \cup \{e^{1,n}_{ij}\}_{0\leq i,j \leq 2}$ generate $e^{0,n-1}_{00}N_n e^{0,n-1}_{00}=e^{1,n-1}_{00}N_n e^{1,n-1}_{00}$, for all $n\geq 1$. 

\vskip.05in 

For $n=1$, this is condition 2$^\circ$ of Corollary 5.3 applied to $N_1$. Assume we have made this construction up to some $n$. By applying 5.3.2$^\circ$ to 
the SSG II$_1$ factor $P=e^{0,n}_{00}N_{n+1} e^{0,n}_{00}=e^{1,n}_{00}N_{n+1} e^{1,n}_{00}$, one gets two $3 \times 3$ matrix units 
$\{e^{0,n+1}_{ij}\}_{0\leq i,j \leq 2}$, $\{e^{1,n+1}_{ij}\}_{0\leq i,j \leq 2}$ in $P$ that have the same diagonal, $e^{0,n+1}_{ii}=e^{1,n+1}_{ii}$, $0\leq i\leq 2$, 
and together generate $P$. 

Now note that $\cup_{n=1}^m (\{e^{0,n}_{ij}\}_{0\leq i,j \leq 2} \cup \{e^{1,n}_{ij}\}_{0\leq i,j \leq 2})$ generate $N_m$. Also, if we denote by $B^m_k$ 
the algebra generated by $\cup_{n=1}^m (\{e^{k,n}_{ij}\}_{0\leq i,j \leq 2}$, $k=0, 1$, then $B^m_0\simeq B^m_1\simeq \Bbb M_3(\Bbb C)^{\otimes m}$. 

Thus, $R_0 =\vee_m B_0^m$ and $R_1=\vee_m B_1^m$ are hyperfinite II$_1$ factors and they generate a von Neumann subalgebra of $M$ that contains 
all $N_m, m\geq 1$, so they generate $M$. Thus, $M$ is generated by two hyperfinite subfactors $R_0, R_1$ satisfying condition $4^\circ$ in Corollary 5.3. 
\hfill $\square$

\heading 6.  SSG, mean-value property, and tightness  
\endheading

For convenience, let us first recall the conjectures (5.1 $(a)$ and $(b)$ in [P18]) that we wanted to comment on in this section. 

\vskip.05in

\noindent
{\bf 6.1. Conjecture}. $(a)$ {\it If a } II$_1$ {\it factor $M$ has the SSG property then it has a properly infinite $R$-bimodule decomposition, i.e., 
it contains a pair of hyperfinite subfactors $R_0, R_1 \subset M$ such that $R_0 \vee R_1^{op}\subset \Cal B(L^2M)$ is a purely infinite 
von Neumann algebra}.  

$(b)$ {\it If a } II$_1$ {\it factor $M$  has the SSG property, then it has an $R$-tight decomposition, i.e., it contains a pair 
of hyperfinite subfactors $R_0, R_1 \subset M$ such that $R_0 \vee R_1^{op}=\Cal B(L^2M)$}.  

\vskip.05in  

While  6.1.$(b)$ implies 6.1.$(a)$, each one of these statements forces different strategies for constructing recursively the pair $R_0, R_1\subset M$, as a limit of finite 
dimensional factors $R_{0,n} \nearrow R_0$, $R_{1,n} \nearrow R_1$.

Thus, Proposition 2.7 shows that in order to prove Conjecture 6.1.$(a)$ one needs to construct  recursively two increasing sequences of finite dimensional 
factors $R_{0,m}$, $R_{1,m}$, $m\geq 1$, inside $M$, so that if we take $\{\xi_k\}_k\subset L^2M$ to be a dense sequence in the set of vectors 
of $\| \ \|_2$-norm equal to $1$, then at each ``next step'' $n$, the finite dimensional factors $R_{0,n}\supset R_{0,n-1}, R_{1,n}\supset R_{1,n-1}$ 
have to be constructed in $M$ so that the restrictions of the vector states $\omega_{\xi_k}$ 
to $R_{0,n}\vee R_{1,n}^{op}=R_{0,n}\otimes R_{1,n}^{op}$ have Radon-Nykodim derivative with respect to $\tilde{\tau}=\tau \otimes \tau$, denoted $a_k$, 
to satisfy that ``most of its mass'' (i.e.,  its $\tilde{\tau}$-trace) is concentrated on a projection of small $\tilde{\tau}$-trace, for all $1\leq k \leq n$. 
This amounts to requiring that $1-\tilde{\tau}(a_ke_{[1, \infty)}(a_k)) \leq \varepsilon_{n}$, $\forall 1\leq k \leq n$, 
for some constants  $\varepsilon_n \searrow 0$ (e.g., $\varepsilon_n=2^{-n}$). 

To have these conditions satisfied it would be sufficient that at each step $n$ the support of $s(\varphi_k) \in R_{0,n}\otimes R_{1,n}^{op}$ 
is majorized by a projection of the form $\sum_i f_i \otimes g_i$ 
where $f_i \in R_{0,n}$, $g_i \in R_{1,n}$ are projections which on the right give a partition of $1$, $\sum_i g_i=1$, while $f_i$ satisfy $\xi_i g_i=f_i \xi_k g_i$ 
(or at least approximately, in an appropriate sense) for all $1\leq k\leq n$ and 
so that the set $J$ of all $i$'s for which one has $\tau(f_i)\leq \varepsilon_n$ satisfies $\sum_{i\in J} \tau(g_i)\geq 1- \varepsilon_n$. 

This implies that the left support $l(\xi_kg_i)$ of $\xi_kg_i$ should have expectation on $R_{0,n}$ that's supported ``in large part'' 
by a projection of trace $\leq \varepsilon_n$, for all $k$ and all $i\in J$. 
Since the partition $\{g_i\}_i$ can be made with arbitrarily small projections 
independently of the set $\xi_1, ..., \xi_n$, the trace of left supports $\vee_k l(\xi_kg_i)$ can indeed be made small. But it seems quite difficult 
to use the SSG assumption to prove the existence of the 
the finite dimensional factors $R_{0,n}, R_{1,n}$ so that these left supports ``avoid coarseness'', i.e., so that $l(\xi_kg_i)$ avoid being orthogonal 
to a finite dimensional factor $R_{1,n}$ that contains $R_{0,n-1}$.

By contrast, the {\it tightness conjecture} 6.1.$(b)$ seems more prompt to a ``dynamical'' approach, as  suggested by the tightness criteria in Section 3.  
Paradoxically, although this is a stronger statement, it seems more feasible and more intuitive. 

Thus, in this case one would like to construct the increasing sequences of finite dimensional factors $R_{0,n}, R_{1,n}\subset M$ 
so that by averaging over the unitaries in $R_{0,n}, R_{1,n}^{op}$ at each step $n$, one obtains that ``larger and larger'' finite subsets of a countable dense 
subset $\{x_k\}_k$ of the space of trace-class operators on $L^2M$ that have trace equal to $0$, 
$L_0^1(\Cal B(L^2M), Tr):=\{x\in \Cal B(L^2M)\mid \|x\|_{1,Tr}=Tr(|x|)< \infty, Tr(x)=0\}$, get ``more and more anihilated''.  Indeed, by (Proposition 2.5 in [P19]) 
this condition is equivalent to the fact that $(\cup_n R_{0,n}\cup \cup_n R_{1,n}^{op})'\cap \Cal B(L^2M)=\Bbb C$, meaning that $R_0=\overline{\cup_n R_{0,n}}^w$, 
$R_1=\overline{\cup_n R_{1,n}}^w$ give an $R$-tight decomposition of $M$. 

Note that $R$-tightness implies that  the  left-right action of the entire unitary group $\Cal U(M)$ on $\Cal B(L^2M)$ satisfies the following mean-value type property, already  
considered in (Section 7 of [P19]): 

\vskip.05in 
\noindent
{\bf 6.2. Definition}. A II$_1$ factor $M$ has the {\it MV-property} if 
the weak closure of the convex hull of $uv^{op}T{v^{op}}^*u^*$, over unitaries $u, v \in M$, intersects the scalars, for any $T \in \Cal B=\Cal B(L^2M)$. 
Equivalently (cf 2.5 in [P19]), if the $\| \ \|_{1,Tr}$-norm closure of the convex hull of $uv^{op}x{v^{op}}^*u^*$, over unitaries $u, v \in M$, contains $0$, 
for any $x\in L^1_0(\Cal B, Tr)$.

\vskip.05in 

Proving that $M$ is $R$-tight amounts to showing that $\Cal U(M)$ has a pair of ``hyperfinite directions'' along which the ergodicity of 
$\Cal U(M) \times \Cal U(M^{op})$ $ \curvearrowright^{\text{\rm Ad}} \Cal B$ is being realized.  
To construct such a pair recursively, it seems essential to first establish that $M$ has the MV-property. This motivated us to formulate in an initial version of this paper 
the following

\vskip.05in

\noindent
{\bf 6.3. Problem} (cf also 7.4 in [P19])  {\it Do all factors have the MV property}? {\it Do free group factors have the MV-property}? {\it Do SSG factors have the MV-property}? 
\vskip.05in

At that time, it wasn't clear at all whether the SSG property was to play a role in answering this question. 
But the above problems were answered a couple of months later by Das and Peterson in [DPe19], as a consequence of their {\it double ergodicity} theorem. 

Let us briefly describe here  their result. Let  $\Cal U=\{u_k\}_{k\geq 1} \subset \Cal U(M)$ be 
a countable self-adjoint set of  unitary elements with $\Cal U''=M$ and $\alpha=(\alpha_k)_k$ positive elements with $\sum_k \alpha_k=1$. 
For  $T\in \Cal B=\Cal B(L^2M)$ denote  $\varphi(T)=\sum_{k\geq 1} \alpha_k u_kTu_k^*$ 
(resp.  $\varphi^{op}(T)= \sum_{k\geq 1} \alpha_k u^{op}_iT{u_i^{op}}^*$) the averaging of $T$ by the unitaries in $\Cal U$ (resp. in $\Cal U^{op}=J\Cal UJ$) with weights $\alpha$. 
The Das-Peterson double ergodicity theorem states that the averaging of any $T\in \Cal B$ by $\phi=(\varphi + \varphi^{op})/2$ and its powers, 
intersects the scalars, i.e., $\lim_{n\rightarrow \omega} \frac{1}{n} \sum_{k=1}^n \phi^k(T)\in \Bbb C1$, $\forall T\in \Cal B$, where $\omega$ is some fixed free ultrafilter on $\Bbb N$.  

So, more than showing that ANY separable (equivalently countably generated) II$_1$ factor has the MV-property, the Das-Peterson theorem provides  some very concrete identification 
of the corresponding left-right averagings.

For instance, if $M=L(\Bbb F_n)$, then one can take $\varphi$ to be the  Markov c.p. map $\varphi(T)=\frac{1}{2n}\sum_{i=1}^n (u_iTu_i^* + u_i^*Tu_i^*)$, 
where $u_1, ..., u_n\in L(\Bbb F_n)$ are the free generators. 

A particular case of interest is when $M$ is assumed to be SSG and we view it as 
 generated by a pair of triadic factors $B_0, B_1 \simeq \Bbb M_3(\Bbb C)$ that have the same 
diagonal, as in Section 5. Let us  denote by $E, E^{op}$ the Tr-preserving expectations from $\Cal B(L^2M)$ 
onto $(B_0\vee B_1^{op})'$ resp. $(B_1\vee B_0^{op})'$ obtained by 
averaging the operators in $\Cal B(L^2M)$ over $B_0 \vee B_1^{op}$, resp. $B_1, B_0^{op}$. Note that $B_0 \vee B_1^{op}$ are $9 \times 9$ matrix algebras 
and $E, E^{op}$ can each be viewed as averaging by a group of $81=9^2$ many unitaries.  

\proclaim{6.4. Corollary (to [DPe19])} With the above notation for $E, E^{op}$ and $\phi=(E+E^{op})/2$, we have that 
$\lim_{n\rightarrow \omega} \frac{1}{n} \sum_{k=1}^n \phi^k(T)\in \Bbb C1$, $\forall T\in \Cal B(L^2M)$. 
\endproclaim

The  concreteness of the averagings in the MV-property, as provided by the Das-Peterson double ergodicity, with the many options and flexibility for taking the averagings,  
could be very useful for  the tightness conjecture. Note also the  interesting fact that, since MV-property holds in general, if the tightness conjecture is to hold true, 
then the SSG property needs to be used only in the ``second part'' of a potential proof.

In such a {\it Step} 2, one should deduce from the MV-property that there exist two finite dimensional factors $Q_0, Q_1$, that refine a  
previously constructed pair of finite dimensional factors, so that the averaging  over $Q_0 \vee Q_1^{op}$  ``diminishes'' the $\| \ \|_{1,Tr}$-norm of a given finite set $F\subset L^1_0(\Cal B, Tr)$ by a fixed 
universal constant $c<1$.  To do so, one has to argue that if 
for some finite set $J$ of pairs of indices, the average $\| \frac{1}{|J|} \sum_{(i,j) \in J} u_iv_j^{op} x u^*_i{v_j^{op}}^*\|_{1,Tr}$, $x\in F$, is small, then there exists a  subset $J_0\subset J$ 
with $\{u_i, v_j \mid (i,j)\in J_0\}$ generating finite dimensional algebras (in other words, there exist  ``finite dimensional directions'') so that 
$\| \frac{1}{|J_0|} \sum_{(i,j) \in J_0} u_iv_j^{op} x u^*_i{v_j^{op}}^*\|_{1,Tr}$ $< c$, $\forall x\in F$. 

Usually, this is being done by first finding an ``abelian direction'' (a finite partition of $1$ by projections) 
on a ``corner'' of $M$ (and of all inclusions involved), like for instance in [P81]. Viewing the terms of the average as 
sitting on a circle, this means that from the average inequality $\| \frac{1}{|J|} \sum_{(i,j)\in J} u_iv_j^{op} x u^*_i{v_j^{op}}^*\|_{1,Tr} < \varepsilon$, 
one  should deduce that there two points on the circle, i.e., $u_{i_0}, u_{i_1}$, $v_{j_0}, v_{j_1}$,  such that $\|\frac{1}{2} (u_{i_0}v^{op}_{j_0}xu_{i_0}^*{v_{j_0}^{op}}^*
+ u_{i_1}v^{op}_{j_1}xu_{i_1}^*{v_{j_1}^{op}}^*)\|_{1,Tr} < c$, $x\in F$. If this were to be the case, then $u_{i_0}u_{i_1}^*$ and $v_{j_0}v_{j_1}^*$ would give 
the two ``abelian directions'' that would allow moving on with the iterative construction of the finite dimensional algebras. 

This is how the iterative construction in [P81] is being done, to obtain an ergodic hyperfinite copy of $R$ in $M$. But unlike [P81], where the ambient norm is 
the $L^2$-norm $\| \ \|_2$ on $M$, in this attempted construction of $R_0, R_1\subset M$ with $R_0\vee R_1^{op}$ ergodic in $\Cal B(L^2M)$  we are dealing 
with the $L^1$ norm $\| \ \|_{1,Tr}$. So the intuition about existence of two points on the circle with its mid-point diminishing the $L^1$-norm 
no longer works. This is the major hurdle for accomplishing this {\it Step} 2.

Finally, {\it Step} 3 would consist in doing all the above recursively, to get the pair of increasing sequences of finite dimensional factors $R_{0,n}, R_{1,n}\subset M$ that satisfy the desired properties, 
somewhat like in the proof of Theorem 5.4. The fact that SSG is a stable property is crucial for this iterative procedure, because at each step $n$ one has to apply again 
{\it Steps} 1 and 2 to a ``small corner'' of $M$.

There is an alternative way to approach the tightness conjecture, where one  
would exploit the SSG property and Das-Peterson double ergodicity more directly, building the increasing sequences $R_{0,n}, R_{1,n}$ by choosing at each step between 
 two triadic factors with same diagonal  that generate $M$ (or rather, a corner of it), using  a spectral gap type property like Lemma 5.1,  in the spirit of  
the proof of  Theorem 5.4. 
This kind of approach would of course need a ``uniform spectral gap'', that would be the same on each corner of $M$. 

The additional assumption that the fundamental group of $M$ is $\Bbb R_+$, or just  $1/3 \in \Cal F(M)$ ($M^{1/3} \simeq M$), may be useful in this respect. Indeed, this would allow 
taking an isomorphism $\theta: M \simeq pMp$ for $\tau(p)=1/3$, choosing two triadic subfactors with same diagonal $B_0, B_1\subset M$ that generate $M$, and then taking 
$\theta^n(B_0)$, $\theta^n(B_1)$ as generators on the corners at steps $n=1, 2, ....$. Note that at each of these steps one would still have the possibility of conjugating $B_0$, $B_1$ by  
random unitaries on both left and right, an operation that may be of additional  help with the ``uniform spectral gap'' requirement. 

Of course, by the same arguments explained above, solving the tightness conjecture under the assumption $1/3 \in \Cal F(M)$ would still imply that $L(\Bbb F_\infty)$ 
is infinitely generated and that all $L(\Bbb F_t), 1< t \leq \infty$, follow non-isomorphic.

\vskip.05in
\noindent
{\bf 6.4. Remark.} Note that if the tightness conjecture $6.1.(b)$ holds true, then it would also follow that if a II$_1$ factor $M$ has a properly infinite $R$-bimodule 
decomposition, then it has an $R$-tight decomposition. So a  natural test for the tightness conjecture is to prove this implication, independently of Conjectures 6.1. 
We have been able to verify it in many cases (such as for all factors that are weakly thin, in the sense of [GP98], in particular factors with Cartan subalgebras, s-MASAs, 
property Gamma factors, etc), but through rather laborious long proofs, which did not seem of much interest for the scope of this paper and will be detailed elsewhere. 

\vskip0.5in

\head  References \endhead

\item{[AP17]} C. Anantharaman, S. Popa: ``An introduction to II$_1$ factors'', \newline www.math.ucla.edu/$\sim$popa/Books/IIun-v13.pdf

\item{[B05]} M. Brin: {\it On the Zappa-Szep product}, Comm. in Algebra {\bf 33} (2005), 393-424. 

\item{[C75]} A. Connes: {\it Classification of injective factors}, Ann. of Math. {\bf 104} (1976), 73-115.

\item{[DPe19]} S. Das, J. Peterson: {\it Poisson boundaries of} II$_1$ {\it factors}, preprint December 2019.

\item{[Dy93]} K. Dykema: {\it Interpolated free group factors}, Pacific J. Math. 163 (1994), 123-135.

\item{[DyR99]} K. Dykema, F R\u adulescu: {\it Compressions of free products of von Neumann algebras.} Math. Ann., {\bf 316} (2000), 61-82.

\item{[DSSW07]} K. Dykema, A. Sinclair, R. Smith, S. White: 
{\it Generators of} II$_1$ {\it factors}, Oper. Matrices {\bf 2} (2008), 555-582. 

\item{[FM77]} J. Feldman, C.C. Moore: {\it Ergodic equivalence
relations, cohomology, and von Neumann algebras} II, Trans. AMS {\bf 234} (1977), 325-359. 

\item{[GeP96]} L. Ge, S. Popa: {\it On some decomposition properties for factors of type} II$_1$,  Duke Math. J., {\bf 94} (1998), 79-101.

\item{[H15]} B. Hayes:  1-{\it bounded entropy and regularity problems in von Neumann algebras},  arXiv:1505.06682, to appear in Int. Math. 
Res. Notices. 

\item{[IPP05]} A. Ioana, J. Peterson, S. Popa:
{\it Amalgamated Free Products of w-Rigid Factors and Calculation of their Symmetry Groups},
Acta Math. {\bf 200} (2008), No. 1, 85-153. (math.OA/0505589)

\item{[IPV10]} A. Ioana, S. Popa, S. Vaes: {\it A Class of superrigid group von Neumann algebras}, Annals of Math. 
{\bf 178} (2013), 231-286 (math.OA/1007.1412).

\item{[J82]} V.F.R. Jones: {\it Index for subfactors}, Invent. Math. {\bf 72} (1983), 1-25.

\item{[J00]} V.F.R. Jones : {\it Ten problems}, in ``Mathematics: perspectives and frontieres'', pp. 79-91, AMS 2000, V. Arnold, M. Atiyah,
P. Lax, B. Mazur Editors.

\item{[K67]} R.V. Kadison: {\it Problems on von Neumann algebras}, Baton Rouge Conference 1967, unpublished manuscript.

\item{[KrV16]} A-S. Krogager, S. Vaes: {\it Thin} II$_1$ {\it factors with no Cartan subalgebras}, to appear in Kyoto J. of Math., arXiv:1611.02138

\item{[MvN36]} F. Murray, J. von Neumann: {\it On rings of operators}, Ann. Math. {\bf 37} (1936), 116-229.

\item{[MvN43]} F. Murray, J. von Neumann: {\it On rings of operators IV}, Ann. Math. {\bf 44} (1943), 716-808.

\item{[OP03]} N. Ozawa, S. Popa: {\it Some prime factorization results for type} II$_1$ {\it factors}, Invent. Math., {\bf 156}
(2004), 223-234 (math.OA/0302240).

\item{[OP07]} N. Ozawa, S. Popa: {\it On a class of} II$_1$ {\it
factors with at most one Cartan subalgebra}, Annals of Mathematics {\bf 172} (2010),
101-137 (math.OA/0706.3623)

\item{[PiP84]} M. Pimsner, S. Popa: {\it Entropy and index for
subfactors}, Annales Scient. Ecole Norm. Sup., {\bf 19} (1986), 57-106.

\item{[PiV81]} M. Pimsner, D. Voiculescu: {\it $K$-groups of reduced crossed products by free groups}, J. Operator Theory {\bf 8} (1982), 131-156. 

\item{[P81]} S. Popa: {\it On a problem of R.V. Kadison on maximal
abelian *-subalgebras in factors}, Invent. Math., {\bf 65} (1981),
269-281.

\item{[P82]} S. Popa: {\it Notes on Cartan subalgebras in type} II$_1$ {\it factors}. Mathematica Scandinavica, {\bf 57} (1985), 171-188.

\item{[P83]} S. Popa: {\it Semiregular maximal abelian *-subalgebras and the solution to the factor state Stone-Weierstrass problem},
Invent. Math., {\bf 76} (1984), 157-161.

\item{[P94]} S. Popa, {\it Symmetric enveloping algebras,
amenability and AFD properties for subfactors}, Math. Res.
Letters, {\bf 1} (1994), 409-425.

\item{[P97]} S. Popa: {\it Some properties of the symmetric enveloping algebras
with applications to amenability and property T},
Documenta Mathematica, {\bf 4} (1999), 665-744.

\item{[P01]} S. Popa: {\it On a class of type} II$_1$ {\it factors with
Betti numbers invariants}, Ann. of Math {\bf 163} (2006), 809-899
(math.OA/0209310; MSRI preprint June 2001).

\item{[P16]}  S. Popa:  {\it Constructing MASAs with prescribed properties}, Kyoto J. of Math,  {\bf 59} (2019), 367-397, math.OA/1610.08945

\item{[P18]} S. Popa: {\it Coarse decomposition of} II$_1$ {\it factors}, math.OA/1811.11016   

\item{[P19]} S. Popa: {\it On ergodic embeddings of factors}, math.OA/1910.06923

\item{[PS18]} S. Popa, D. Shlyakhtenko: {\it Representing interpolated free group factors as group factors}, math.OA/1805.10707, to appear in 
Groups, Geormetry, and Dynamics.

\item{[PV11]} S. Popa, S. Vaes: {\it Unique Cartan decomposition for} II$_1$ 
{\it factors arising from arbitrary actions of free groups}, Acta Mathematica, {\bf 194} (2014), 237-284. 

\item{[Po67]} R. Powers: {\it Representations of uniformly hyperfinite algebras and their associated von Neumann rings}, Ann. Math. {\bf 86} (1967), 138-171. 

\item{[R91]} F. Radulescu: {\it The fundamental group of the von Neumann algebra of a free group with infinitely many generators is $\Bbb R_+$}, J. Amer. Math. Soc. {\bf 5} (1992), 517-532. 

\item{[R92]} F. Radulescu: {\it Random matrices, amalgamated free products and subfactors of the von Neumann algebra of a free group, of noninteger 
index}, Invent. Math. {\bf 115} (1994), 347-389.
  
\item{[Sa71]} H. Sakai: ``C$^*$-algebras and W$^*$-algebras'', Springer-Verlag, Berlin-Heidelberg-New York, 1971.

\item{[Sh05]} J. Shen: {\it Type} II$_1$ {\it factors with a single generator} , J. Oper. Theory {\bf 62} (2009), 421-439.   

\item{[S97]} D. Shlyakhtenko: {\it Some applications of freeness with amalgamation}, J. reine angew. Math. {\bf 500} (1998), 191-212

\item{[V88]} D. Voiculescu: {\it Circular and semicircular systems and free product factors}, Prog. in Math. {\bf 92}, Birkhauser, Boston, 1990, pp. 45-60.

\item{[V96]} D. Voiculescu: {\it The analogues of entropy and Fisher's information measure in free probability theory} III: {\it absence of Cartan subalgebras}, 
GAFA {\bf 6} (1996), 172-199. 

\enddocument